\documentclass[11pt]{amsart}

\usepackage{amsmath,amssymb,amsthm}
\usepackage[margin=1in]{geometry}
\usepackage{microtype}      
\usepackage{hyperref,array}
\hypersetup{colorlinks=true,linkcolor=blue,citecolor=blue,urlcolor=blue}

\theoremstyle{plain}
\newtheorem{theorem}{Theorem}[section]

\newtheorem{thmx}{Theorem}                  

\newtheorem{corollary}[theorem]{Corollary}
\newtheorem{proposition}[theorem]{Proposition}
\newtheorem{lemma}[theorem]{Lemma}
\theoremstyle{definition}

\newtheorem*{definition*}{Definition}
\newtheorem{example}[theorem]{Example}
\theoremstyle{remark}
\newtheorem{remark}[theorem]{Remark}

\DeclareMathOperator{\Fix}{Fix}

\DeclareMathOperator{\Isom}{Isom}
\DeclareMathOperator{\Min}{Min}

\DeclareMathOperator{\gdim}{dim_{\!g}}
\newcommand{\R}{\mathbb{R}}
\newcommand{\N}{\mathbb{N}}
\newcommand{\CAT}{\ensuremath{\operatorname{CAT}}}
\newcommand{\UEG}{\mathrm{UEG}}
\newcommand{\bd}{\partial_\infty}

\title[Uniform elliptic gaps and fixed points]
{A local-to-global fixed point theorem for $\CAT(0)$ spaces with a uniform elliptic gap
}
\author{Jeroen Schillewaert}
\address{Department of Mathematics\\
University of Auckland\\
38 Princes Street\\
1010, Auckland, New Zealand}
\email{j.schillewaert@auckland.ac.nz}
\date{\today}

\begin{document}

\begin{abstract}
A finitely generated group acting locally elliptically on a complete
\CAT(0) space of finite geometric dimension fixes a point if distances
to elliptic fixed sets are uniformly controlled by displacements. For arbitrary
groups, the same hypothesis gives a fixed point in the space or its
visual boundary. The proof combines a first-moment random-walk
displacement estimate with Izeki's zero-drift theorem.
For finite-rank $\R$-buildings whose Euclidean metric factor has
finite linear image, we obtain an all-type local-to-global theorem,
including the nondiscrete $\widetilde G_2$ case, with an explicit gap
and without metric completeness; the linear-image restriction can be
removed from the fixed-point conclusion when the Euclidean factor has
dimension at most three.
A tangent-cone criterion characterizes the gap and yields a
cubical bound with optimal asymptotic growth $d/\pi$ in dimension $d$.
For Hadamard manifolds, the gap is equivalent to uniformly bounded
finite stabilizers in the isometry image, with stabilizer bounds
of optimal polynomial growth in the gap constant.
\end{abstract}

\maketitle

\section{Introduction}\label{sec:intro}

If $G$ acts by isometries on a metric space $X$ and every element of $G$ fixes a point of $X$, that is, the action is \emph{locally elliptic}, under which hypotheses on $G$ and $X$ does $G$ have a global fixed point? Our hypothesis is a uniform displacement estimate comparing the
distance to an elliptic fixed-point set with the displacement of
the corresponding isometry. Actions need not be faithful; statements about orders or cardinalities of isometries refer to
the image of the action. We write $gx$ for the image of $x$ under $g$. For $g\in G$, set  $\Fix_X(g)=\{x\in X:gx=x\}$ and let $\Fix_X(G)=\bigcap_{g\in G}\Fix_X(g)$.
An element $g$ is \emph{elliptic} if $\Fix_X(g)\ne\varnothing$. Let $C>0$.
The action satisfies $\UEG(C)$ if
\begin{equation}\label{eq:ueg}
 d(x,\Fix_X(g))\leq C\,d(x,gx)
\end{equation}
for every elliptic $g\in G$ and every $x\in X$, where
$d(x,D)=\inf_{y\in D}d(x,y)$ for a nonempty subset $D$ of $X$.
The action has a \emph{uniform elliptic gap} if it satisfies $\UEG(C)$ for
some $C>0$. For a single map on a Hilbert space, \eqref{eq:ueg} is the standard linear-regularity inequality in the terminology of
Bauschke--Noll--Phan \cite[Definition 2.1]{BNP15}.
Here the constant is uniform over
all elliptic elements of the isometry image.
This condition depends only on the isometry image.
If $\UEG(C)$ holds with $C<1/2$, the image contains no nonidentity
elliptic element, since $d(x,gx)\leq2d(x,\Fix_X(g))$.
For an \(\R\)-tree, one can take \(C=1/2\). The condition allows
elliptic elements of infinite order: the binary adding-machine
automorphism of the rooted binary tree fixes the root and has
order \(2^n\) on level \(n\), hence has infinite order.

For locally elliptic actions the uniform gap implies a local-to-global fixed point theorem. 
\begin{thmx}[Local-to-global fixed points]
\label{thm:main}
Let \(X\) be a complete \CAT(0) space of finite geometric
dimension. Let a group \(G\) act isometrically and locally
elliptically on \(X\) with a uniform elliptic gap.
If \(G\) is finitely generated, then \(G\) fixes a point of \(X\).
For arbitrary \(G\), there is a fixed point in \(X\cup\bd X\).
\end{thmx}

For bounded-exponent groups, the fixed-point conclusion is due to
Hou--Izeki--Ji--Karlsson--Wu \cite[Theorem 1.2(2)]{HIJKW26}.
Theorem~\ref{thm:main} replaces the order bound by a metric condition
on the given action, allowing infinite-order elliptic isometries.
The gap bounds the distances from the basepoints to the individual fixed sets
in rescaled ultralimits and forces the locally elliptic Euclidean image arising in the proof to be finite, see Lemma \ref{lem:euclidean}.
For the point-or-boundary conclusion, finite geometric dimension
can be replaced by properness or, for countable groups, by finite
telescopic dimension (Corollary~\ref{cor:point-or-boundary}).

The uniform elliptic gap hypothesis cannot be omitted from the general
local-to-global statement.
Bass constructed a finitely generated locally elliptic subgroup
of \(\Isom(\R^4)\) with no common fixed point
\cite[Counterexample 1.10, pp.~19--20]{Bas80}; see also Breuillard and Fujiwara
\cite[Example 9.1]{BF21}. Finite generation is needed for a fixed
point in \(X\) even for complete trees
\cite[II.7, Example 7.11]{BH99}. Nor can the dimension restriction
be dropped: Osajda constructs unbounded cubical actions of certain
finitely generated Burnside groups \cite[Theorems 1 and 2]{Osa18}.
The actions extend isometrically to the metric completion, where
their orbits remain unbounded and bounded exponent yields the
gap by Proposition~\ref{prop:order}. Moreover, finite telescopic
dimension does not suffice for the point conclusion: the same
groups admit unbounded locally elliptic actions with a uniform
gap on a complete \CAT(-1) space of telescopic dimension one
(Example~\ref{ex:telescopic}).
Conversely, a fixed point does not imply a uniform elliptic gap:
the cyclic group generated by an irrational rotation of $\R^2$
fixes the origin, while a sequence of its nonidentity powers
approaches the identity.

The tree case is classical: Serre's theorem for simplicial trees
\cite[I.6.5]{Ser80} and the extension to real trees by Morgan and Shalen
\cite[Proposition II.2.15]{MS84} give the finitely generated
local-to-global conclusion.
Norin, Osajda and Przytycki prove the analogous conclusion for finitely
generated groups acting by metric-preserving simplicial automorphisms
on complete \CAT(0) triangle complexes with finitely many isometry
classes of triangles incident to each vertex, provided that either elliptic
elements have finite order, or the complex is locally finite, or the
angles are rational with respect to the action in their sense
\cite[Theorem 1.1]{NOP22}.

For buildings we use Parreau's displacement argument
\cite[Theorem 4.1 and Section 4.2]{Par00}, also developed by
Ramos-Cuevas \cite[Proposition 2.1]{RC15}. Applied to an
elliptic fixed set, the apartment-projection argument gives
the gap without metric completeness. We include a direct
proof with an explicit constant depending only on the reduced
spherical Coxeter complex (Proposition~\ref{prop:building-gap}).
The same conclusion holds with a Euclidean factor when its
linear image is finite. Combining this with our local-to-global
theorem and Struyve's fixed-point theorem for bounded actions
gives the following result in all types.

\begin{thmx}[$\mathbb{R}$-buildings]\label{thm:buildings}
Let $X=X_0\times\R^k$ be a finite-rank $\R$-building, with
$\R^k$ its Euclidean metric factor. Let a finitely generated
group $G$ act isometrically on $X$, and suppose its induced
linear image on $\R^k$ is finite. If every element fixes a point
of $X$, then $G$ fixes a point of $X$.
Metric completeness of $X$ is not required.
\end{thmx}

The finite-linear-image hypothesis is automatic for actions
preserving a finite essential Euclidean Weyl structure: the
reflecting normals span the Euclidean factor, and the linear parts
preserve the corresponding finite reflection arrangement.
It cannot be omitted in arbitrary dimension, as Bass's example
shows. Lemma~\ref{lem:euclidean-obstruction} identifies the Euclidean
factor as the only obstruction and removes this hypothesis when
that factor has dimension at most three.

Theorem~\ref{thm:buildings} is new for nondiscrete
$\widetilde G_2$-buildings and recovers the theorem of the author, Struyve and Thomas
\cite[Theorem A]{SST25} for $\R$-buildings of types $\widetilde A_2$ and $\widetilde C_2$.
In higher rank, it extends the finitely generated fixed-point
conclusion of Parreau \cite[Theorem 1 and Corollary 3]{Par03}
to metric building actions without a realization of the acting
group in a valued-field linear group; it does not retain
Parreau's more general hypothesis of generation by a bounded set.
For discrete buildings, it recovers the theorem of Haettel and
Osajda for the classical families
$\widetilde A_n,\widetilde B_n,\widetilde C_n,\widetilde D_n$
\cite[Corollary B(1)]{HO25}. The triangle-complex theorem
\cite[Theorem 1.1]{NOP22} covers discrete rank-two affine buildings,
including type $\widetilde G_2$.

An action on a \CAT(0) space $X$ is called \emph{non-elementary}
if it has no bounded orbit in $X$ and no finite orbit in $\bd X$,
as in \cite[p.~2]{LLS24}. If a non-elementary action satisfies the finite-linear-image
hypothesis of Theorem~\ref{thm:buildings}, then $k=0$.
Indeed, for $k>0$ the embedded sphere
$\partial_\infty\R^k\subseteq\partial_\infty X$ has finite
orbits: translations act trivially on it, and the linear
image is finite.
In $\R$-buildings of type
$\widetilde A_2$, Le Bars--L\'ecureux--Schillewaert show that, for
non-elementary type-preserving actions and symmetric generating
measures with finite second moment, strongly regular hyperbolic
elements occur with probability tending to one. Their argument
also gives the local-to-global theorem \cite[Theorems C and D]{LLS24},
without metric completeness. For finite-rank $\R$-buildings
without Euclidean metric factors, our random-walk consequence is
axiality without metric completeness (Corollary~\ref{cor:axial});
strong regularity in general type requires further assumptions,
as Quint's example shows (Remark~\ref{rem:quint}).

\medskip

\paragraph{Further applications.}
We prove using Izeki's drift theorem \cite[Theorem A]{Ize23} that on a complete proper or finite-telescopic-dimensional $\mathrm{CAT}(0)$ space a uniform elliptic gap and the absence of fixed points in $X\cup \partial_\infty X$ or invariant flats imply that symmetric generating walks with finite second displacement moment have positive drift and ellipticity probability tending to zero (Theorem~\ref{thm:random}).
The gap also holds for actions by cubical
automorphisms of finite-dimensional \CAT(0) cube
complexes: in Section~\ref{sec:cubes} we prove a linear
displacement bound with optimal asymptotic coefficient $1/\pi$
in the dimension. The proof uses the same-constant tangent-cone
criterion of Theorem~\ref{thm:tangent-gap}, which characterizes
the gap through the tangent cones of elliptic fixed sets.
The resulting fixed-point conclusion is already known from
Sageev's work \cite[proof of Theorem 5.1]{Sag95} and Leder--Varghese
\cite[Theorem A]{LV20}; see also Fioravanti's finite-rank median
generalization \cite[Theorem 3.1]{Fio18} and the extension by
Genevois, Lonjou and Urech \cite[Theorem 1.3]{GLU24}. Section~\ref{sec:examples} discusses bounded elliptic
orders, bounded point stabilizers, geometric actions,
and permanence under products and finite-index overgroups.

For geometric actions of nonamenable discrete groups,
symmetric generating walks with finite first displacement
moment have positive drift and axial random elements with
probability tending to one, without a dimension assumption
(Corollary~\ref{cor:geometric}). This applies in particular
to non-elementary geometric actions.

On a Hadamard manifold,
the gap is equivalent to point stabilizers in the isometry
image being finite with uniformly bounded cardinalities
(Proposition~\ref{prop:smooth}). This gives a quantitative
finite-image conclusion for locally elliptic actions,
without finite generation (Corollary~\ref{cor:smooth}). 

In this setting the gap is equivalent to a uniform lower bound
on the nonzero singular values of the derivatives of elliptic
elements. Jordan's theorem then gives stabilizer bounds with
optimal growth exponent in the gap constant for fixed
manifold dimension
(Proposition~\ref{prop:jordan} and
Remark~\ref{rem:smooth-sharpness}).
Table~\ref{tab:constants} collects the gap constants.

\paragraph{Method of proof.}
The probabilistic argument extends the block interchange strategy developed in
\cite[Lemma 3.2]{HIJKW26}. Under a finite first displacement
moment, without symmetry, we obtain
\(d(o,Z_n^2o)/n\to2\ell\) in \(L^1\), where \(\ell\) is the drift.
The scalar argument uses the gap only through the single-isometry
estimate of Lemma~\ref{lem:ratio-gap}: a uniform strict upper bound below \(2\)
for \(d(o,g^2o)/d(o,go)\) when \(g\) is elliptic and moves \(o\).
Thus positive drift forces the probability of ellipticity to tend
to zero. With symmetry and finite second moment, Izeki's
zero-drift theorem \cite[Theorem A]{Ize23} implies positive drift
under the geometric hypotheses of Theorem~\ref{thm:random}.

\smallskip

For the fixed-point theorem we combine the normalization argument going back to Shalom \cite[Lemma 6.3]{Sha00}, in the
form given by Bourdon \cite[Lemma 3.2]{Bou12}, in the $\mathrm{CAT}(0)$ form used
by Izeki–Karlsson \cite[Proposition 8]{IK26}, with the boundary dimension-reduction argument of Hou–Izeki–Ji–
Karlsson–Wu \cite[Lemmas 2.1--2.2, Theorem 3.7 and Proposition 3.8]{HIJKW26}.
The gap preserves individual fixed points in moving-basepoint
ultralimits and forces a locally elliptic Euclidean image to be
finite, allowing the argument to apply to the given action without
a group-wide Euclidean fixed-point assumption. The usual hypothesis $\mathrm{Hom}(G,\R) = 0$ rules out translations in the line factor of the
boundary reduction; here local ellipticity makes that translation character vanish directly.

\smallskip

For buildings, we use the displacement argument of Parreau \cite[p.38-39]{Par00}, also developed by Ramos-
Cuevas \cite[Proposition 2.1]{RC15}. Applied to the fixed set of an elliptic isometry, the apartment-projection argument
gives the gap without metric completeness. We include a direct proof with an explicit
constant depending only on the reduced spherical Coxeter complex. A finite linear image yields the gap on the Euclidean metric factor. Semisimplicity in that
generality is due to Le Bars--L\'ecureux--Schillewaert
\cite[Theorem E]{LLS24}, extending Parreau's theorem
\cite[Th\'eor\`eme 4.1 and Corollaire 4.2]{Par00}.
Struyve's enlargement and fixed-point theorems
\cite[Lemma 4.4 and Main Result 1]{Str11} then yield a fixed point
in the original, possibly incomplete, building.

We conclude with an overview of the gap constants we obtain in this paper.
\begin{table}[ht]
\centering
\renewcommand{\arraystretch}{1.25}
\begin{tabular}{
  >{\raggedright\arraybackslash}p{0.42\textwidth}
  >{\raggedright\arraybackslash}p{0.27\textwidth}
  >{\raggedright\arraybackslash}p{0.20\textwidth}}
\hline
Action & Gap constant & Reference \\
\hline
$\R$-trees & $1/2$ & Introduction \\
Complete $\CAT(0)$ spaces, elliptic orders at most $N\geq 2$
  & $(2\sin(\pi/N))^{-1}$ & Proposition~\ref{prop:order} \\
Cubical automorphisms, dimension $d$, Euclidean metric
  & $C_d$ & Proposition~\ref{prop:cube-intro} \\
Cubical automorphisms, dimension $d$, piecewise $\ell^1$ metric
  & $d/2$ & Proposition~\ref{prop:cube-intro} \\
Median graphs of cube dimension $d$, fixed-vertex automorphisms
  & $d/2$ & Theorem~\ref{thm:median-gap} \\
Reduced $\R$-buildings
  & $C_0$ from \eqref{eq:building-constant}
  & Proposition~\ref{prop:building-gap} \\
Finite products of actions with constants $C_i$
  & $\max_i C_i$ & Section~\ref{sec:examples} \\
Index-$m$ overgroups, $m\geq 2$, complete $\CAT(0)$ spaces
  & $mC+(2\sin(\pi/m))^{-1}$
  & Proposition~\ref{prop:finite-index} \\
\end{tabular}
\medskip
\caption{Elliptic gap constants. The Euclidean cubical
constant $C_d$ is defined in \eqref{eq:cubical-constant}
and has optimal asymptotic growth $d/\pi$.
The median-graph and piecewise $\ell^1$ constants are sharp.
The building constant is explicit but coarse.}
\label{tab:constants}
\end{table}

\section{Background}\label{sec:background}

\subsection{Convexity, projections, and elliptic isometries}

Let $X$ be a $\CAT(0)$ space, not necessarily complete, and let
$D\subseteq X$ be a nonempty complete convex subspace. Then every point of $X$
has a unique nearest point in $D$, and the resulting projection
$p_D\colon X\to D$ is $1$-Lipschitz \cite[Prop.~II.2.4]{BH99}. By uniqueness,
$p_D$ commutes with every isometry of $X$ that preserves $D$. This applies in
particular to nonempty closed convex subsets of a complete $\CAT(0)$ space,
and to flats, which are complete and convex. For an isometry $g$ of $X$, the
fixed-point set $\operatorname{Fix}(g)$ is always closed and convex
\cite[Prop.~II.2.2]{BH99}, and it is complete whenever $X$ is complete.

Every nonempty bounded subset $B$ of a complete $\CAT(0)$ space has a unique
circumcenter, that is, the center of the unique closed ball of minimal radius
containing $B$ \cite[Prop.~II.2.7]{BH99}. By uniqueness, the circumcenter is
fixed by every isometry preserving $B$. In particular, we have the
\emph{Bruhat--Tits fixed point theorem}: if a group acts by isometries on a
complete $\CAT(0)$ space with a bounded orbit, then it has a global fixed
point; see
\cite[Cor.~II.2.8]{BH99}.

For an isometry $g$ of a \CAT(0) space $(X,d)$, put
\[
 |g|=\inf_{x\in X}d(x,gx),\qquad
 \Min(g)=\{x\in X:d(x,gx)=|g|\}.
\]

An isometry $g$ of $X$ is \emph{semisimple} if its displacement function
$x\mapsto d(x,gx)$ attains its infimum, i.e. \(\Min(g)\neq\varnothing\). A semisimple isometry is either
elliptic, i.e.\ it has a fixed point, or \emph{axial} (also called hyperbolic):
there are an isometric embedding $a\colon\mathbb{R}\to X$ and a number $t_g>0$
such that $g\,a(t)=a(t+t_g)$ for all $t\in\mathbb{R}$. An isometry that is not
semisimple is called \emph{parabolic} \cite[II.6]{BH99}.

\subsection{Boundary points and parallel sets}
All rays are parametrized by arc length. The \emph{visual boundary}
$\bd X$ is the set of asymptotic classes of rays in $X$. If $X$ is complete, each
boundary point has a unique representative ray starting at each
point of $X$. For an isometry $g$ of a \CAT(0) space $X$ fixing
$\xi\in\bd X$ and a ray $c$ representing $\xi$, convexity gives
\begin{equation}\label{eq:raybound}
 d(c(t),gc(t))\leq d(c(0),gc(0))\quad(t\geq0).
\end{equation}

The \emph{parallel set} of a geodesic line in a \CAT(0) space is the
union of the lines at finite Hausdorff distance from it. If $X$ is a
complete \CAT(0) space, the parallel set of any line in $X$ is closed
and convex and splits isometrically as $P=Y\times\R$, where $Y$
is a complete \CAT(0) space. Its vertical lines are exactly the lines
parallel to the original line.
If a group $G$ acts isometrically on such an $X$, preserves this
parallel class, and fixes its positive endpoint, its action on $P$
has the form
\begin{equation}\label{eq:parallel-action}
 g(y,t)=(\bar g y,t+\beta(g)),
\end{equation}
where $g\mapsto\bar g$ is an isometric action on $Y$ and
$\beta:G\to(\R,+)$ is a homomorphism.
See \cite[II.2, Theorem 2.14, and II.8]{BH99}. The definitions of the visual boundary and of a parallel set do
not require completeness. Completeness is used here to obtain
representative rays from arbitrary basepoints and the closed
parallel-set decomposition with complete transverse factor.

\subsection{Dimension and ultralimits}\label{sec:dim-ultra}

For a complete $\CAT(0)$ space $X$, let $\gdim X$ denote
Kleiner's geometric dimension. Kleiner proved that
\[
\gdim X=\sup\{\dim_{\rm cov}K:K\subseteq X
\text{ is compact}\},
\]
where $\dim_{\rm cov}$ is Lebesgue covering dimension
\cite[Theorem A(1)]{Kle99}. In particular, a complete convex
subspace has no larger geometric dimension. We shall also use
\[
\gdim(Y\times\mathbb R)=\gdim Y+1
\]
for nonempty complete $\CAT(0)$ spaces $Y$; this is the usual
product formula for geometric dimension (see the proof of
Lemma A.12 in the erratum appended to \cite{CL10}). A uniform
finite upper bound on geometric dimension is preserved under
pointed ultralimits \cite[Proposition 6]{IK26}; see also
\cite[Section 2.1]{CL10}. Rescaling the metric does not change
geometric dimension. A nonempty zero-dimensional $\CAT(0)$
space is a point, since a segment between distinct points
would be a compact subset of covering dimension one.
A complete $\mathbb R$-building has geometric dimension equal
to its apartment dimension: Kleiner and Leeb 
bound the relative homology of every open pair above that
dimension \cite[Corollary 6.1.1, p.~172]{KL97}, and Kleiner identifies the resulting
bound with geometric dimension; an apartment gives equality
\cite[Theorem A(2)]{Kle99}.
A $\CAT(0)$ cube complex of bounded cube dimension $d$, with
unit Euclidean cubes, is complete and has geometric dimension
$d$. Completeness follows from the finite-shape theorem
\cite[I.7, Theorem 7.50]{BH99}; the local cone description
\cite[I.7, Theorems 7.16 and 7.39]{BH99}, applied
also to links, bounds the successive spaces of directions by
the cube dimension, and a $d$-cube gives the reverse inequality.

For a nonprincipal ultrafilter $\omega$ on $\N$, write
$X_\omega=\lim_\omega(X_i,d_i,p_i)$ for the pointed ultralimit:
its points are classes of \emph{admissible} sequences, those with
$\sup_i d_i(x_i,p_i)<\infty$, identified at zero distance, and
$d_\omega([x_i],[y_i])=\lim_\omega d_i(x_i,y_i)$.
It is complete, even if the factors are not
\cite[I.5.53]{BH99}, and is \CAT(0) when the factors are
\cite[II.3.10(2)]{BH99}; the dimension bound above applies to
complete \CAT(0) factors. If $G$ acts on the factors and
$\sup_i d_i(p_i,gp_i)<\infty$ for each $g\in G$, then
$g[x_i]=[gx_i]$ defines an isometric action on $X_\omega$.
For finitely generated $G$, it suffices to check this bound on
a finite generating set. Rescaled metrics $d_i=\lambda_i d$,
with $\lambda_i>0$, are allowed; admissibility and displacement
bounds are measured using $d_i$.

\subsection{Zero drift}
Let a countable group \(G\) act isometrically on a \CAT(0) space
\((X,d)\), let \(o\in X\), and let \(\mu\) be a probability measure
on \(G\).
On the probability space $(G^\N,\mu^{\otimes\N})$, let $s_i$ be
the coordinate maps and write $Z_n=s_1\cdots s_n$. Thus the
increments are independent with law $\mu$; probabilities are denoted
by $\mathbb P$.
For \(p=1,2\), finite \(p\)-th moment means
\(\sum_{g\in G}\mu(g)d(o,go)^p<\infty\); the condition is independent
of \(o\).  Under the finite first moment assumption,
Kingman's subadditive ergodic theorem \cite{Kin68} gives a deterministic drift \(\ell\geq0\)
such that
\begin{equation}\label{eq:drift}
 \frac{d(o,Z_no)}n\longrightarrow\ell
 \quad\text{almost surely and in }L^1.
\end{equation}
It does not depend on \(o\).  The measure is symmetric if
\(\mu(g)=\mu(g^{-1})\), and generating if its support generates \(G\).
For a finitely supported measure all moments of displacement are finite.
If $G$ is generated by a finite symmetric set $S$, then
$d(o,go)\leq M|g|_S$, where $M=\max_{s\in S}d(o,so)$.
Thus a finite word-length moment implies the corresponding
displacement moment.

For an integer $n\geq0$, a \CAT(0) space $(X,d)$ has
\emph{telescopic dimension} at most $n$ if every asymptotic cone
has geometric dimension at most $n$. An asymptotic cone of $(X,d)$
is an ultralimit $\lim_\omega(X,\lambda_i d,p_i)$, where $\omega$
is a nonprincipal ultrafilter on $\N$, $p_i\in X$, and
$\lambda_i>0$ with $\lambda_i\to0$. 
In the terminology of Caprace and Lytchak, finite geometric
dimension implies finite telescopic dimension \cite[Section 2.1]{CL10}.

We use the following formulation of a structural alternative proved by Izeki \cite[Theorem A]{Ize23}.
\begin{theorem}[Izeki]\label{thm:izeki}
Let \(X\) be a complete \CAT(0) space which is proper or has
finite telescopic dimension. Let a countable group \(G\) act
isometrically on \(X\), and let \(\mu\) be a symmetric generating
probability measure on $G$ with finite second displacement moment.
If the drift is zero, then \(G\) fixes a point of \(X \cup \bd X\) or
preserves a Euclidean flat in \(X\).
\end{theorem}

\subsection{Tangent cones and the elliptic gap}
For $p\in X$, let $\Sigma_pX$ be the space of directions at $p$ with the
angle metric, and let $T_pX$ be the completion of the Euclidean cone over
$\Sigma_pX$, with vertex $0$ and metric $d_p$; write $\log_p x$ for the vector
of length $d(p,x)$ in the direction of $[p,x]$.
The map $\log_p$ is $1$-Lipschitz by angle comparison.
For a convex set $F$ containing $p$, let $T_pF$ be the closure
in $T_pX$ of the cone generated by $\log_p F$.
An isometry $g$ fixing $p$ induces a homogeneous isometry
$D_pg$ of $T_pX$. See \cite[II.3.18--II.3.19]{BH99}.

\begin{theorem}[Tangent-cone criterion]\label{thm:tangent-gap}
Let a group act isometrically on a complete \CAT(0) space $X$,
let $C>0$ and for each elliptic $g$, let $F=\Fix_X(g)$. The action satisfies $\UEG(C)$ if and only if,
for every elliptic $g$, every $p\in F$, and every
$v\in T_pX$,
\begin{equation}\label{eq:tangent-gap}
 d_p(v,T_pF)\leq C\,d_p(v,D_pg\,v).
\end{equation}
In particular, the gap implies
$T_p\Fix_X(g)=\Fix_{T_pX}(D_pg)$.
\end{theorem}
\begin{proof}
Suppose first that $\UEG(C)$ holds. For $v=\log_p x$, let
$x_t$ be the point at parameter $t\in[0,1]$ on $[p,x]$.
For every $q\in F$, angle comparison gives
\[
 d_p(v,T_pF)\leq t^{-1}d_p(\log_p x_t,\log_p q)
 \leq t^{-1}d(x_t,q).
\]
Taking the infimum over $q$, applying the gap, and letting
$t\downarrow0$ yields
\[
 d_p(v,T_pF)\leq C\lim_{t\downarrow0}
             \frac{d(x_t,gx_t)}t
 =C\,d_p(v,D_pg\,v).
\]
Homogeneity and density give the inequality on all of $T_pX$.
The inclusion $T_pF\subseteq\Fix(D_pg)$ is automatic; the reverse
inclusion follows from the inequality and closedness of $T_pF$.

Conversely, let $p$ be the projection of $x$ onto $F$ and put
$v=\log_p x$. The projection criterion
\cite[II.2, Proposition 2.4]{BH99} says that the angle from $v$
to every nonzero vector of $T_pF$ is at least $\pi/2$.
Consequently $d_p(v,T_pF)=d(p,x)$, and
\[
 d(x,F)=d_p(v,T_pF)
 \leq C\,d_p(v,D_pg\,v)\leq C\,d(x,gx).
\]
\end{proof}

\begin{remark}
The forward implication does not require completeness. For the
converse it is enough that, for every elliptic element $g$ and
$x\in X$, the nonempty set $\Fix_X(g)$ contains a nearest point
to $x$. In particular, the conclusion holds if every nonempty
elliptic fixed set is complete, even when $X$ itself is incomplete.
\end{remark}

The criterion uses the tangent cone of the actual fixed set.
Replacing it by $\Fix(D_pg)$ requires that fixed tangent vectors
be approximable by directions of actual fixed segments.
This condition is automatic in the following setting. We call $X$ \emph{locally conical} if, for every $p\in X$, there
is $r_p>0$ such that
\[
\log_p:B_X(p,r_p)\longrightarrow B_{T_pX}(0,r_p)
\]
is an isometry onto the indicated ball.

\begin{corollary}[Locally conical spaces]\label{cor:conical-gap}
Suppose that $X$ is complete, \CAT(0), and locally conical.
An isometric action on $X$ satisfies $\UEG(C)$ if and only if
the action of each point stabilizer on $T_pX$ satisfies
$\UEG(C)$, with the same constant for all $p$.
This applies, in particular, to piecewise Euclidean \CAT(0)
polyhedral complexes with finitely many isometry types of cells.
\end{corollary}
\begin{proof}
The local isometry is equivariant for the stabilizer of $p$.
Every $D_pg$-fixed vector therefore gives an actual fixed
segment, so $T_p\Fix_X(g)=\Fix_{T_pX}(D_pg)$.
Apply Theorem~\ref{thm:tangent-gap}.
For the polyhedral assertion, apply
\cite[I.7, Theorem 7.16]{BH99} after barycentric subdivision;
see also \cite[I.7, Proposition 7.49 and Theorem 7.50]{BH99}.
\end{proof}

General $\R$-buildings need not be locally conical;
see Remark~\ref{rem:buildings-polyhedral}.

\section{The scalar random-walk argument}\label{sec:scalar}

The block interchange below is the one used in the bounded argument of \cite[Lemma 3.2]{HIJKW26}.  The following deterministic consequence of the $\mathrm{CAT}(0)$ midpoint inequality is the geometric
input for the first-moment argument.

\begin{lemma}[Block interchange]\label{lem:two-block}
Let $(X,d)$ be a \CAT(0) space, let $o\in X$, and let
$A,B\in\Isom(X)$. Put $L(h)=d(o,ho)$ for $h\in\Isom(X)$, and set
\[
 a=L(A),\quad b=L(B),\quad s=L(AB),\qquad
 E=\frac{a^2+b^2}{2}-\frac{s^2}{4}.
\]
Then \(E\geq0\) and
\begin{equation}\label{eq:two-block}
 L((AB)^2)\geq 2L(BA)-4\sqrt E.
\end{equation}
\end{lemma}
\begin{proof}
Let \(p=Ao\), \(q=B^{-1}o\), \(z=ABo\), \(w=(AB)^{-1}o\).
In the Euclidean comparison triangles for $(o,p,z)$ and
$(o,q,w)$, the quantity $E$ is the squared length of the
median to the side of length $s$.
For the midpoints \(m\in[o,z]\), \(m'\in[o,w]\), we obtain
\(d(p,m)^2,d(q,m')^2\leq E\).
Also \(d(m,m')\leq d(z,w)/2\).  Hence
\[
 L(BA)=d(p,q)\leq2\sqrt E+\tfrac12d(z,w)
                 =2\sqrt E+\tfrac12L((AB)^2).
\]
\end{proof}

\begin{proposition}\label{prop:ratio}
Let a countable group $G$ act isometrically on a \CAT(0) space $(X,d)$,
let $o\in X$, and put $L(g)=d(o,go)$. Let $\mu$ be a probability
measure on $G$ such that $\sum_g\mu(g)L(g)<\infty$.
For independent increments $s_i$ of law $\mu$, set
$Z_n=s_1\cdots s_n$ and let
$\ell=\lim_{n\to\infty}n^{-1}L(Z_n)$ be the almost sure drift. Then
\[
 \frac{L(Z_n^2)}n\longrightarrow2\ell
 \quad\text{in $L^1$}.
\]
If \(\ell>0\), then \(L(Z_n^2)/L(Z_n)\to2\) in $L^p$ for every $1\leq p < \infty$,
where the ratio is set to zero when \(L(Z_n)=0\). The latter event occurs only finitely often almost surely.
\end{proposition}
\begin{proof}
For $n\geq2$, write $k=\lfloor n/2\rfloor$,
$A_n=s_1\cdots s_k$, and $B_n=s_{k+1}\cdots s_n$.

Kingman's theorem gives $L(Z_n)/n\to\ell$ in $L^1$ and
\[
 \left(\frac{L(A_n)}n,\frac{L(B_n)}n,
       \frac{L(A_nB_n)}n,\frac{L(B_nA_n)}n\right)
 \longrightarrow
 \left(\frac\ell2,\frac\ell2,\ell,\ell\right)
 \quad\text{in probability}.
\]
Indeed, the first two coordinates have the laws of walks of
lengths $k$ and $n-k$, while both products have the law of $Z_n$
by invariance of the increment law under cyclic permutation.

Let $E_n$ denote the quantity $E$ in
Lemma~\ref{lem:two-block} for $A_n,B_n$. The preceding convergences and continuity of squaring give
$\frac{E_n}{n^2}\longrightarrow 0$
in probability.
By Lemma~\ref{lem:two-block} $L(Z_n^2)\geq 2L(B_nA_n)-4\sqrt{E_n}$, thus $\frac{L(Z_n^2)}n\longrightarrow2\ell$
in probability, using also $L(Z_n^2)\leq2L(Z_n)$.

Since $L(Z_n)/n\to\ell$ in $L^1$, the family
$\{2L(Z_n)/n:n\geq1\}$ is uniformly integrable. The inequality
\[
0\leq\frac{L(Z_n^2)}n\leq\frac{2L(Z_n)}n
\]
therefore implies that $\{L(Z_n^2)/n:n\geq1\}$ is uniformly
integrable. Hence the convergence upgrades to convergence in $L^1$.
If $\ell>0$, division by $L(Z_n)/n\to\ell$ in probability
gives
\[
\frac{L(Z_n^2)}{L(Z_n)}\longrightarrow2
\qquad\text{in probability}.
\]

Put $R_n=L(Z_n^2)/L(Z_n)$, with $R_n=0$ if $L(Z_n)=0$. By
subadditivity, $L(Z_n^2)\leq 2L(Z_n)$, so $0\leq R_n\leq 2$ and
$|R_n-2|\leq 2$. Hence, for $1\leq p<\infty$ and $\varepsilon>0$,
\[
 \mathbb E\,|R_n-2|^p\leq \varepsilon^p+2^p\,\mathbb P\bigl(|R_n-2|>\varepsilon\bigr).
\]
Letting $n\to\infty$ and then $\varepsilon\downarrow0$, convergence in
probability gives $R_n\to2$ in $L^p$.

\end{proof}

\begin{lemma}[Elliptic displacement gap]\label{lem:ratio-gap}
Let $(X,d)$ be a \CAT(0) space, let $g\in\Isom(X)$ be elliptic,
let $x\in X$, and let $C\geq1/2$. If
$d(x,\Fix_X(g))\leq C d(x,gx)$, then
\begin{equation}\label{eq:ratio-gap}
 d(g^{-1}x,gx)\leq2q_Cd(x,gx),\qquad
 q_C=
 \begin{cases}
 C,&1/2\leq C\leq1/\sqrt2,\\
 \sqrt{1-1/(4C^2)},&C\geq1/\sqrt2.
 \end{cases}
\end{equation}
In particular, $q_C<1$.
\end{lemma}
\begin{proof}
Put $L=d(x,gx)$ and $D=d(g^{-1}x,gx)$; assume $L>0$.
For $\varepsilon>0$, choose $p\in\Fix_X(g)$ with
$R=d(x,p)\leq(C+\varepsilon)L$, and let $m$ be the midpoint of
$[g^{-1}x,gx]$. \CAT(0) comparison gives
\[
 d(x,m)\leq\sqrt{L^2-D^2/4},\qquad
 d(p,m)\leq\sqrt{R^2-D^2/4},
\]
and hence $R-\sqrt{R^2-D^2/4}\leq\sqrt{L^2-D^2/4}$.
Put $c=C+\varepsilon$ and $t=D/(2L)$.
Since $D\leq2R$ and
$r\mapsto r-\sqrt{r^2-D^2/4}$ is nonincreasing for $r\geq D/2$,
\[
 t\leq c,\qquad c-\sqrt{c^2-t^2}\leq\sqrt{1-t^2}.
\]
If $c\leq1/\sqrt2$, the first inequality gives $t\leq c$.
If $c\geq1/\sqrt2$ and $c\geq\sqrt{1-t^2}$, squaring the
second inequality gives $t^2\leq1-1/(4c^2)$.
In the remaining case $c<1$, and
$t^2<1-c^2\leq1-1/(4c^2)$, since $c^4\geq1/4$.
Letting $\varepsilon\downarrow0$ completes the proof.
\end{proof}

\begin{remark}[Sharpness]
The bound is sharp, even for finite cyclic actions, for every $C \ge 1/2$.
Put $\theta = 2\arcsin\bigl(1/(2C)\bigr)$ and choose $m \ge 4$ with $m\theta \ge 2\pi$.
Rotation by $\theta$ on the $\mathrm{CAT}(0)$ Euclidean cone of total angle $m\theta$ \cite[II.3, Theorem 3.14]{BH99}
generates a cyclic group of order $m$ with optimal gap constant $C$.
The condition $m \ge 4$ makes the angular separation of $g^{-1}x$ and $gx$
equal to $2\theta$ before truncation at $\pi$. Thus, for $d(x,0)=R$,
\[
  d(x,gx) = 2R\sin\frac{\theta}{2},
  \qquad
  d(g^{-1}x,gx) = 2R\sin\frac{\min\{\pi,\,2\theta\}}{2},
\]
and
\[
\frac{d(g^{-1}x,gx)}{d(x,gx)}
=\frac{\sin(\min\{\pi,2\theta\}/2)}{\sin(\theta/2)}
=2q_C.
\]
In particular, $\mathbb{R}^2$ gives equality at $C = (2\sin(\pi/m))^{-1}$
for $m \ge 4$.
\end{remark}

\begin{corollary}\label{cor:zero-drift}
Let a countable group \(G\) act isometrically on a \CAT(0) space \(X\),
and suppose that the action satisfies a uniform elliptic gap.
Let \(\mu\) be a probability measure on \(G\) with finite first moment,
and let \(Z_n\) be the random walk with independent increments of law
\(\mu\).  If its drift \(\ell\) is positive, then
\[
 \mathbb P(Z_n\text{ is elliptic})\longrightarrow0.
\]
Consequently, if the action is locally elliptic, then \(\ell=0\).
If every element of \(G\) acts semisimply and \(\ell>0\), then
\(\mathbb P(Z_n\text{ is axial})\to1\).
\end{corollary}
\begin{proof}
Choose $o\in X$, put $L(g)=d(o,go)$, and choose a gap constant
$C\geq1/2$. Suppose $\ell>0$, and let $R_n=L(Z_n^2)/L(Z_n)$,
with $R_n=0$ when $L(Z_n)=0$. Lemma~\ref{lem:ratio-gap} gives
\[
 \mathbb P(Z_n\text{ is elliptic})
 \leq\mathbb P(R_n\leq2q_C)\longrightarrow0,
\]
by Proposition~\ref{prop:ratio}, since $2q_C<2$. The other
assertions follow immediately.
\end{proof}

The numerical constants below are chosen for convenience, we did not optimise them.

\begin{corollary}[Transfer of displacement concentration]\label{cor:concentration}
In the setting of Corollary~\ref{cor:zero-drift}, suppose that
the gap constant is $C\geq1/2$ and the drift $\ell$ is positive.
Put $\delta=1-q_C$ and $\eta=\ell\delta^2/32$.
If constants $A,b>0$ satisfy
\[
 \mathbb P\bigl(|d(o,Z_jo)-j\ell|>\eta j\bigr)\leq Ae^{-bj}
 \qquad(j\geq1),
\]
then, for every integer $n\geq16/\delta^2$,
\[
 \mathbb P(Z_n\text{ is elliptic})\leq
       4A e^{-b\lfloor n/2\rfloor}.
\]
If $\mu$ is supported on $S$ and the right-hand side is less than
one, some word of length at most $n$ in $S$ is nonelliptic.
\end{corollary}
\begin{proof}
Use the two blocks $A_n,B_n$ of Proposition~\ref{prop:ratio},
and put $\epsilon=2\eta=\ell\delta^2/16$.
Outside four exceptional events, for walks of lengths
$\lfloor n/2\rfloor,\lceil n/2\rceil,n,n$, respectively, we have
\[
 \frac{L(A_n)}n,\frac{L(B_n)}n\leq\frac\ell2+\epsilon,
 \qquad
 \ell-\epsilon\leq\frac{L(A_nB_n)}n,
                      \frac{L(B_nA_n)}n\leq\ell+\epsilon.
\]
The parity error in the first two bounds is at most
$\ell/(2n)\leq\eta$. Lemma~\ref{lem:two-block} gives
$E_n/n^2\leq3\ell\epsilon/2+3\epsilon^2/4\leq2\ell\epsilon$,
and therefore
\[
 \frac{L(Z_n^2)}n
 \geq2(\ell-\epsilon)-4\sqrt{2\ell\epsilon}
 >2q_C(\ell+\epsilon)\geq2q_C\frac{L(Z_n)}n.
\]
Lemma~\ref{lem:ratio-gap} excludes ellipticity on this event.
A union bound gives the asserted probability estimate.
\end{proof}

\begin{theorem}[Random walks and the elliptic gap]\label{thm:random}
Let $X$ be a complete \CAT(0) space which is proper or has finite
telescopic dimension. Let a countable group $G$ act isometrically on
$X$ with a uniform elliptic gap. Suppose that $G$ fixes no point of
$X \cup \bd X$ and preserves no Euclidean flat in $X$.
Let $\mu$ be a symmetric probability measure on $G$ whose support
generates $G$, and suppose that
$\sum_{g\in G}\mu(g)d(o,go)^2<\infty$ for a basepoint $o\in X$.
For independent increments $s_i$ of law $\mu$, set
$Z_n=s_1\cdots s_n$. The almost sure drift
$\ell=\lim_{n\to\infty}n^{-1}d(o,Z_no)$ is positive, and
\[
 \mathbb P(Z_n\text{ is elliptic})\longrightarrow0.
\]
\end{theorem}
\begin{proof}[Proof of Theorem~\ref{thm:random}]
The hypotheses exclude both alternatives in
Theorem~\ref{thm:izeki}, so the drift is positive.
Corollary~\ref{cor:zero-drift} now gives the conclusion.
\end{proof}

\section{The local-to-global fixed-point theorem}\label{sec:proof}

We combine the Shalom--Bourdon normalization argument
\cite[Lemma 6.3]{Sha00}\cite[Lemma 3.2 and Proposition 3.1]{Bou12},
in the \CAT(0) form used by Izeki and Karlsson
\cite[Proposition 8]{IK26}, with the boundary
dimension-reduction argument of Hou, Izeki, Ji, Karlsson and Wu
\cite[Lemmas 2.1--2.2 and Proposition 3.8]{HIJKW26}.
Theorem~\ref{thm:izeki} leaves a fixed boundary point or an invariant
flat when the drift is zero. To handle these alternatives for the given action, we must control
individual fixed points when rescaling the metric, restricting to
invariant complete convex subspaces, and taking ultralimits with
moving basepoints. The gap gives this control for the
normalization ultralimit; for the ultralimit along a fixed boundary
ray, fixed rays give the required control directly.
Ellipticity need not pass to a moving-basepoint ultralimit without a uniform bound on the distances from the basepoints to the fixed-point sets.

\paragraph{Why the basepoints matter.}
For $n\geq3$, let the generator of $\mathbb Z$ act on
$X_n=\R^2$ by rotation $R_{2\pi/n}$ about $0$, and take
$p_n=(n/(2\pi),0)$. Every factor action has finite image and is
locally elliptic, while for every fixed $k\in\mathbb Z$,
\[
d(p_n,R_{2\pi k/n}p_n)
=\frac n\pi\left|\sin\frac{\pi k}{n}\right|
\leq |k|.
\]
Thus the action on the pointed ultralimit is defined.
In coordinates centred at $p_n$, its generator is
\[
z\longmapsto R_{2\pi/n}z+(R_{2\pi/n}-I)p_n,
\]
which converges to the translation $z\mapsto z+(0,1)$.
The factor fixed points are not admissible:
$d(p_n,\Fix(R_{2\pi/n}))=n/(2\pi)\to\infty$.
The gap constants in these finite cyclic actions tend to
infinity. A constant uniform in the factors prevents this
loss of ellipticity in Lemma~\ref{lem:permanence}.

\begin{lemma}[Permanence]\label{lem:permanence}
Let $C>0$.
\begin{enumerate}
\item Let a group $G$ act isometrically on a \CAT(0) space $(X,d)$
with $\UEG(C)$. The same gap holds after any positive rescaling
of $d$. It also holds on every nonempty invariant complete convex
subspace $D\subseteq X$. Local ellipticity on $X$ implies local
ellipticity on $D$.
\item Let $\omega$ be a nonprincipal ultrafilter on $\N$. Let a
group $G$ act isometrically and locally elliptically on pointed
\CAT(0) spaces $(X_i,d_i,p_i)$ with $\UEG(C)$. Assume that every
$g\in G$ has uniformly bounded basepoint displacement:
$\sup_i d_i(p_i,gp_i)<\infty$.
The induced action on $\lim_\omega(X_i,d_i,p_i)$ is locally
elliptic and satisfies $\UEG(C)$.
\end{enumerate}
\end{lemma}
\begin{proof}
Rescaling multiplies both sides of the gap inequality by the same
number. For the restriction, let $x\in D$ and let $g$ be elliptic
on $X$. For $\varepsilon>0$ choose $q\in\Fix_X(g)$ with
$d(x,q)\leq C d(x,gx)+\varepsilon$.
Projection onto the complete convex set $D$ is nonexpanding and
commutes with $g$, so $p_Dq\in\Fix_D(g)$ and
$d(x,p_Dq)\leq C d(x,gx)+\varepsilon$.
Letting $\varepsilon\downarrow0$ proves the gap and also proves
preservation of ellipticity under restriction. For the ultralimit, fix $g\in G$ and an admissible representative
$(x_i)$ of $x$. Choose $q_i\in\Fix_{X_i}(g)$ with
\[
 d_i(x_i,q_i)\leq C d_i(x_i,gx_i)+1/i.
\]

Put
\[R=\sup_i d_i(p_i,x_i)<\infty,\qquad
M_g=\sup_i d_i(p_i,gp_i)<\infty.
\]
Then
\[d_i(x_i,gx_i)\leq2R+M_g,\qquad
d_i(p_i,q_i)\leq R+C(2R+M_g)+1.
\]
Thus $(q_i)$ is admissible.
Its class $q=[q_i]$ is fixed by $g$, and
$d_\omega(x,q)\leq C d_\omega(x,gx)$.
This applies to every $g$ and every $x$.
\end{proof}

The local ellipticity hypothesis is essential here: for general actions, new elliptic elements
in the limit need not satisfy the same gap.

A classical theorem of Jordan states that, for every $n\geq 1$, there is a constant $J(n)$ such
that every finite subgroup of $\mathrm{GL}_n(\mathbb C)$ contains an abelian normal subgroup of index at
most $J(n)$. We refer to Breuillard's recent modern exposition \cite[Theorem 1.1]{Bre23}.\footnote{Collins
\cite{Col07} proves $J(n)=(n+1)!$ for $n\geq71$.
Embedding in $\mathrm{GL}_{71}(\mathbb C)$ also gives the
fully explicit bound $J(n)\leq(\max\{n,71\}+1)!$ in all dimensions;
smaller-dimensional optimal values are given in that reference.}

For $C\geq1/2$, put
\[
 Q(C)=\left\lfloor\frac{\pi}{\arcsin(1/(2C))}\right\rfloor,
 \qquad
 \mathcal B_n(C)=
 \min\left\{(1+2C)^{n^2},
 J(n)\max_{0\leq r\leq\lfloor n/2\rfloor}2^{n-2r}Q(C)^r\right\};
\]
for $0<C<1/2$, put $\mathcal B_n(C)=1$.

\smallskip

The integer $Q(C)$ bounds the order of each nontrivial planar
rotation occurring in the abelian subgroup used below.
The quantity $\mathcal B_n(C)$ bounds the size of an orthogonal
group with normal gap constant $C$: its two terms come from
a separation estimate and from Jordan's theorem, respectively.
It will later give uniform bounds on point stabilizers of
Hadamard-manifold actions.

\begin{lemma}[Finite orthogonal groups with a normal gap]\label{lem:orthogonal}
Let $n\geq1$, $C>0$, and $H\leq O(n)$. Suppose that
\[
 \|(I-A)v\|\geq C^{-1}\|v\|
 \qquad(A\in H,\ v\perp\ker(I-A)).
\]
Then $H$ is finite and $|H|\leq\mathcal B_n(C)$.
For $C<1/2$, the group $H$ is trivial.
\end{lemma}
\begin{proof}
By assumption $\|I-A\|_{\mathrm{op}}\geq C^{-1}$ if $A\neq I$. Moreover $\|I-A\|_{\mathrm{op}}\leq 1+\|A\|_{\mathrm{op}}=2$ since $A$ is
orthogonal; thus if $C<1/2$ then $H=I$.
For $A\neq B\in H$ we have
$\|A-B\|_{\mathrm{op}}=\|I-A^{-1}B\|_{\mathrm{op}}\geq C^{-1}$, since $A^{-1}B\in H\setminus\{I\}$ and the operator norm is invariant under
left multiplication by orthogonal matrices.
Open operator-norm balls of radius $r=(2C)^{-1}$ about a finite
subset $F\subseteq H$ are disjoint and lie in the ball of radius
$1+r$ about zero in the $n^2$-dimensional matrix space.
Comparison of the Lebesgue measures of these balls and of the ball of
radius $1+r$ about $0$ in $M_n(\mathbb R)\cong\mathbb R^{n^2}$ gives
$|F|\leq(1+2C)^{n^2}$.
Thus $H$ is finite with the same bound.

Choose an abelian normal subgroup $A\leq H$ of index at most
$J(n)$. Its real orthogonal representation splits into $s$
one-dimensional sign representations and $r$ two-dimensional
rotation representations, with $s+2r=n$.
If a rotation image has order $m>1$, some element acts on that
plane through angle $2\pi/m$. The normal gap gives
$2\sin(\pi/m)\geq C^{-1}$, hence $m\leq Q(C)$.
A trivial rotation image also has order at most $Q(C)$.
The representation is faithful, so
$|A|\leq2^sQ(C)^r$. Multiplying by $[H:A]\leq J(n)$ proves
the polynomial bound.
\end{proof}

\begin{lemma}[Euclidean rigidity, with a cardinality bound]\label{lem:euclidean}
Let $m\geq1$ and $C>0$. If a group acts isometrically and locally
elliptically on $\R^m$ with $\UEG(C)$, its isometry image
$\Gamma$ is finite, satisfies $|\Gamma|\leq\mathcal B_m(C)$,
and fixes a point.
\end{lemma}
\begin{proof}
Write $\gamma x=A_\gamma x+b_\gamma$ and choose
$q\in\Fix(\gamma)$. Then
$\gamma(q+v)=q+A_\gamma v$ and
$\Fix(\gamma)=q+\ker(I-A_\gamma)$.
Applying the gap to $q+v$, for $v\perp\ker(I-A_\gamma)$, gives
$\|(I-A_\gamma)v\|\geq C^{-1}\|v\|$.
Lemma~\ref{lem:orthogonal} bounds the linear image by
$\mathcal B_m(C)$. The kernel of the linear-part map consists of
elliptic translations, hence is trivial.
Thus $\Gamma$ has the same cardinality as its linear image.
The average of a finite $\Gamma$-orbit is a fixed point.
\end{proof}

For a group $G$ with a nonempty finite symmetric generating set
$S$ and an isometric action on a metric space $(X,d)$, define
\[
 f_S(x)=\max_{s\in S}d(x,sx),\qquad
 \delta_S(X)=\inf_{x\in X}f_S(x).
\]
The identity is allowed in $S$, including when $G$ is trivial. Such an action has \emph{almost fixed points} if \(\delta_S(X)=0\), a condition independent of the finite generating set.

We recall Shalom's selection argument \cite[Lemma 6.3]{Sha00}, in the
form given by Bourdon \cite[Lemma 3.2]{Bou12} and used for \CAT(0)
spaces by Izeki--Karlsson \cite[Lemma 7 and Proposition 8]{IK26}.

\begin{lemma}[Normalization]\label{lem:normalization}
Let $X$ be a complete \CAT(0) space.
Let $G$ act isometrically on $X$ without a fixed point, and let
$S$ be a finite symmetric generating set for $G$. Given a
nonprincipal ultrafilter $\omega$ on $\N$, there are $p_i\in X$
and $\lambda_i>0$ such that the action on
\[
 Y=\lim_\omega(X,\lambda_i d,p_i)
\]
is defined, $Y$ is complete \CAT(0), if $\gdim X\leq n<\infty$, then $\gdim Y\leq n$, and
$\delta_S(Y)\geq1/2$.
If the original action is locally elliptic and satisfies
$\UEG(C)$, so does the action on $Y$.
\end{lemma}

\begin{proof}
Put $f(x)=\max_{s\in S}d(x,sx)$. Since
\[
|d(x,sx)-d(y,sy)|\leq2d(x,y),
\]
the function $f$ is $2$-Lipschitz. As $S$ generates $G$ and
the action has no global fixed point, $f(x)>0$ for every $x\in X$.
For every positive integer $i$, some $p_i\in X$ satisfies
\[
 d(p_i,y)\leq i f(p_i)\quad\Longrightarrow\quad
 f(y)\geq\tfrac12 f(p_i).
\]
Otherwise, for this $i$, start with $x_0\in X$ and choose
recursively
\[
 d(x_j,x_{j+1})\leq i f(x_j),\qquad
 f(x_{j+1})<\tfrac12 f(x_j).
\]
Then $f(x_j)\leq2^{-j}f(x_0)$ and
\[
 \sum_{j=0}^\infty d(x_j,x_{j+1})\leq2i f(x_0)<\infty.
\]
Completeness gives a limit $x_\infty$, and continuity gives
$f(x_\infty)=0$, a contradiction.

Choose $p_i$ as above and set $\lambda_i=f(p_i)^{-1}$. In
the metric $d_i=\lambda_i d$,
$\max_{s\in S}d_i(p_i,sp_i)=1$. If $g=s_1\cdots s_k$ is a
word in $S$, the triangle inequality gives
$d_i(p_i,gp_i)\leq k$. Hence the action is defined on the
ultralimit.

If $y=[y_i]\in Y$, admissibility gives $R$ such that
$d_i(p_i,y_i)\leq R$ for all $i$. For $i\geq R$, the choice
of $p_i$ gives $\max_{s\in S}d_i(y_i,sy_i)\geq1/2$.

For finitely many bounded real sequences $(a_{i,s})_i$, we have
\[
\lim_\omega\max_{s\in S}a_{i,s}
=\max_{s\in S}\lim_\omega a_{i,s}.
\]
Indeed, the finitely many defining $\omega$-large sets may be
intersected, or equivalently one may use that the maximum
map on $\mathbb R^S$ is $1$-Lipschitz for the supremum norm.
Applying this to $a_{i,s}=d_i(y_i,sy_i)$ gives
\[
 \max_{s\in S}d_Y(y,sy)
     =\lim_\omega\max_{s\in S}d_i(y_i,sy_i)\geq\tfrac12.
\]
The remaining metric and dimension assertions follow from
ultralimit stability. The additional inheritance assertion is
Lemma~\ref{lem:permanence}(2), applied after rescaling.
\end{proof}

\begin{lemma}[Boundary reduction for locally elliptic actions]\label{lem:boundary}
Let $G$ act isometrically and locally elliptically on a
complete \CAT(0) space $X$. Let $S$ be a nonempty finite symmetric generating
set with $\delta_S(X)>0$. If $G$ fixes $\xi\in\bd X$, there
is a locally elliptic action on a complete \CAT(0)
space $Y$ such that $ \delta_S(Y)\geq\delta_S(X)$.
If $X$ has finite geometric dimension, then
$\gdim Y\leq\gdim X-1$.
If the original action satisfies $\UEG(C)$, the new action
satisfies the same inequality.
\end{lemma}

\begin{proof}
Choose a nonprincipal ultrafilter $\omega$ on $\N$.
Fix $o\in X$ and let $c$ be the unit-speed ray from $o$ to $\xi$.
The rays $c$ and $gc$ are asymptotic, and convexity gives
\[
 d(c(t),gc(t))\leq d(o,go)\qquad(t\geq0).
\]
Put $p_i=c(i)$ and form the pointed ultralimit
\[
Z=\lim_\omega(X,d,p_i).
\]
Thus the metric is unchanged and the basepoint moves to infinity
along $c$. The bound $d(p_i,gp_i)\leq d(o,go)$ shows that
$g[z_i]=[gz_i]$ defines an isometric action on $Z$: if $(z_i)$
is admissible, then
\[
d(gz_i,p_i)\leq d(z_i,p_i)+d(gp_i,p_i)
\]
is uniformly bounded, and isometries preserve the equivalence
relation defining the ultralimit.

For $z=[z_i]\in Z$, finiteness of $S$ gives
\[
 \max_{s\in S}d_Z(z,sz)
  =\lim_\omega\max_{s\in S}d(z_i,sz_i)
  \geq\delta_S(X).
\]

Fix $g\in G$ and choose $q_g\in\Fix_X(g)$. Let $r_g$ be
the unit-speed ray from $q_g$ to $\xi$. Uniqueness of the
based ray implies $gr_g=r_g$, and
\[
 d(c(t),r_g(t))\leq d(o,q_g)\qquad(t\geq0).
\]
Therefore $(r_g(i))$ is an admissible sequence of $g$-fixed
points. Its class is a fixed point in $Z$. 

For $t\in\R$, define $a(t)=[c(i+t)]$, replacing the
terms with $i<-t$ by $c(0)$. These are admissible
sequences and $d_Z(a(t),a(u))=|t-u|$, so $a$ is a line.
For every $g\in G$,
\[
 d_Z(a(t),ga(t))\leq d(o,go)\qquad(t\in\R).
\]
Hence $ga$ is parallel to $a$ with the same orientation.
The parallel set $P(a)$ is a nonempty closed convex invariant
subspace of $Z$, and splits as $Y\times\R$. The action is
\[
 g(y,t)=(\bar g y,t+\beta(g)).
\]

Projection of a $g$-fixed point of $Z$ onto $P(a)$ is
$g$-fixed. Hence every $g$ is elliptic on $P(a)$, and the
displayed product action forces $\beta(g)=0$. The slice
$Y\times\{0\}$ is invariant. A fixed point $(y,t)$ of $g$
on $P(a)$ yields the fixed point $(y,0)$ on the slice.
Its joint displacement is its joint displacement in $Z$,
so the lower bound $\delta_S(X)$ persists.

Finally $Y$ is complete \CAT(0), and if $\gdim X$ is finite then
\[
 \gdim Y+1=\gdim P(a)\leq\gdim Z\leq\gdim X.
\]
If the original action satisfies $\UEG(C)$, Lemma~\ref{lem:permanence}(2)
transfers the gap to $Z$ and Lemma~\ref{lem:permanence}(1) transfers it to
$P(a)$ and the slice. 
\end{proof}

\begin{proof}[Proof of Theorem~\ref{thm:main}]
First suppose that $G$ is finitely generated. Choose a nonempty
finite symmetric generating set $S$, a gap constant $C>0$, and a
nonprincipal ultrafilter $\omega$ on $\N$.
Suppose, for a contradiction, that \(G\) has no fixed point in \(X\).
Lemma~\ref{lem:normalization} gives a locally elliptic
\(\UEG(C)\) action on a complete \CAT(0) space of no larger
geometric dimension, with positive joint minimal displacement.

Whenever the resulting action fixes a boundary point, apply
Lemma~\ref{lem:boundary}. The dimension decreases and the
joint minimal displacement stays positive. Since a nonempty
zero-dimensional \CAT(0) space is a point, this process stops
with an action on a complete \CAT(0) space $Y$ of finite geometric dimension
having no fixed point in \(Y\cup\bd Y\).

There is no invariant flat in \(Y\). Indeed, its restricted action
would be locally elliptic with the same gap by
Lemma~\ref{lem:permanence}, and Lemma~\ref{lem:euclidean}
would give a global fixed point on any positive-dimensional
such flat. A zero-dimensional invariant flat is already a
global fixed point. Both contradict positive joint displacement.
The finitely generated group $G$ is countable. Choose the uniform
probability measure on the nonempty finite
symmetric generating set $S$. Its drift on $Y$ is zero, since the
action is locally elliptic and has a uniform elliptic gap, by
Corollary~\ref{cor:zero-drift}. The space $Y$ is complete and has
finite geometric, hence finite telescopic, dimension. Theorem~\ref{thm:izeki} therefore gives a
$G$-fixed point of $\bd Y$ or a $G$-invariant Euclidean flat in $Y$.
Both possibilities have been excluded. This contradiction proves
the finitely generated assertion.

Now let \(G\) be arbitrary. For each finitely generated subgroup
\(H\leq G\) put \(D_H=\Fix_X(H)\). The assertion just proved makes
\(D_H\) a nonempty closed convex subset of \(X\). These sets form a
filtering family, since
\[
 D_{\langle H_1,H_2\rangle}\subseteq D_{H_1}\cap D_{H_2}.
\]

By the filtering-family theorem of Caprace--Lytchak
\cite[Theorem 1.1]{CL10}, either their intersection is
nonempty or \(\bigcap_H\bd D_H\) is nonempty. The first
intersection is \(\Fix_X(G)\). In the second case, every element
\(g\in G\) fixes every point of \(\bd D_{\langle g\rangle}\),
because it fixes each ray in \(D_{\langle g\rangle}\)
pointwise. Thus any point of \(\bigcap_H\bd D_H\) is \(G\)-fixed.
\end{proof}

\begin{corollary}[Point-or-boundary conclusion]
\label{cor:point-or-boundary}
Let a group $G$ act isometrically and locally elliptically, with
a uniform elliptic gap, on a complete \CAT(0) space $X$.
Suppose either that $X$ is proper, or that $G$ is countable
and $X$ has finite telescopic dimension.
Then $G$ fixes a point of $X\cup\bd X$.
\end{corollary}
\begin{proof}
First suppose that $G$ is countable. Choose a symmetric generating
measure with finite second displacement moment. 

Fix $o\in X$ and choose a sequence $(g_i)_{i\geq1}$ containing
every element of $G$. Put
\[
a_i=\frac{2^{-i}}{1+d(o,g_io)^2},\qquad A=\sum_{i\geq1}a_i,
\qquad
\mu=\frac1{2A}\sum_{i\geq1}a_i(\delta_{g_i}+\delta_{g_i^{-1}}).
\]
Then $0<A\leq1$, the measure is symmetric and has support $G$,
and, since $d(o,g_i^{-1}o)=d(o,g_io)$,
\[
\sum_{g\in G}\mu(g)d(o,go)^2
=\frac1A\sum_{i\geq1}a_i d(o,g_io)^2\leq\frac1A<\infty.
\]
Corollary~\ref{cor:zero-drift} gives zero drift.
Theorem~\ref{thm:izeki} supplies a fixed point in $X\cup\bd X$
or an invariant Euclidean flat. A zero-dimensional invariant flat
is a fixed point. On a positive-dimensional invariant flat,
projection preserves local ellipticity and the gap, and
Lemma~\ref{lem:euclidean} gives a global fixed point.

Now suppose that $X$ is proper and $G$ is arbitrary.
The space $K=X\cup\bd X$, with the cone topology, is compact
Hausdorff, and isometries extend to homeomorphisms of $K$
\cite[II.8.5--II.8.9, pp.~263--264]{BH99}.
For every finitely generated subgroup $H\leq G$, the countable
case gives a nonempty closed set $E_H=\Fix_K(H)$.
These sets have the finite intersection property, since
$E_{\langle H_1,\ldots,H_r\rangle}\subseteq
E_{H_1}\cap\cdots\cap E_{H_r}$.
Compactness gives $\Fix_K(G)=\bigcap_H E_H\ne\varnothing$.
\end{proof}

\begin{example}[The telescopic-dimension obstruction]\label{ex:telescopic}
There is a finitely generated bounded-exponent group acting
locally elliptically with a uniform gap on a complete \CAT(-1)
space of telescopic dimension one, but without a fixed point
in the space.

Indeed, Osajda's construction \cite[Theorems 1 and 2]{Osa18}
gives a finitely generated group $G$ of exponent $N<\infty$
with an unbounded cubical action. The standard halfspace-cocycle
construction of Niblo and Reeves \cite[Section 3, Technical Lemma, pp.~5--6]{NR97}
gives an affine isometric action $g\cdot v=\pi_gv+b_g$ on a
Hilbert space $H$, with $\|b_g\|^2=2d_1(o,go)$ for a chosen
cubical vertex $o$; hence its orbit is unbounded.

Extend this action to the complete hyperbolic Hilbert space
$X=H\times(0,\infty)$ by
$g(v,t)=(\pi_gv+b_g,t)$, using the half-space metric
\[
 \cosh d((v,t),(w,s))
 =1+\frac{\|v-w\|^2+(t-s)^2}{2ts}.
\]
This is the standard horospherical extension of an affine
Hilbert action described by Monod and Py \cite[Section 2.4 and Corollary 4.3]{MP19}.
The space is \CAT(-1), and the distance formula shows that
the induced action is isometric and unbounded.
Every element has order dividing $N$, where $N\geq2$ since
the action is unbounded. Finite-order isometries of a complete
$\CAT(0)$ space are elliptic, and
Proposition~\ref{prop:order} gives
\[
\UEG\bigl((2\sin(\pi/N))^{-1}\bigr).
\]

Finally, $X$ is unbounded and Gromov hyperbolic, and therefore
has telescopic dimension one \cite[Introduction and Section 2.1]{CL10}.
Thus the boundary alternative in
Corollary~\ref{cor:point-or-boundary} cannot in general be removed,
even for finitely generated groups.
\end{example}

\section{\texorpdfstring{$\R$}{R}-buildings}\label{sec:buildings}

An $\R$-building is understood in the sense of Parreau
\cite[Sections 1--2]{Par00}, with its maximal apartment system.
Its apartments are copies of a finite-dimensional Euclidean
space equipped with a finite reflection group. Their
common dimension $r$ is called the rank. No crystallographic or discreteness
assumption is imposed. Thus our conventions include the
nondiscrete rank-two buildings of Berenstein and Kapovich
\cite{BK}. Such a building is \CAT(0), but need not be
metrically complete \cite[Proposition 2.10]{Par00}.
The asymptotic classes of Weyl sectors form a spherical building
at infinity \cite[Section 1.5]{Par00}; its spherical metric is
the angular metric on the visual boundary. The cone topology
on that boundary is a different, also useful, topology.
We distinguish a prescribed Weyl subdivision from the canonical
reduced structure. A panel is thin if it belongs to exactly two
chambers, and thick if it belongs to at least three. Removing
thin subdivisions gives the thick frame of the spherical
building. Schwer and Weniger's reduction theorem
\cite[Theorem B and Section 5]{SW19} gives a metric product
\[
X=X_0\times\R^k,\qquad
\partial_\infty X=\partial_\infty X_0*S^{k-1},
\]
where $\R^k$ is the maximal Euclidean factor and
$\partial_\infty X_0$ has its canonical thick structure.
The construction does not require metric completeness.
The reduced structure, including its spherical factor, is
intrinsic to the spherical metric; see Kleiner and Leeb
\cite[Section 3.3, Proposition 3.7.1 and Corollary 3.7.4]{KL97}.
Hence metric isometries preserve both boundary factors.
In the displayed product, the union of the rays from $(x,z)$
with endpoints in $S^{k-1}$, together with $(x,z)$, is
$\{x\}\times\R^k$, and the
union of the rays with endpoints in $\partial_\infty X_0$,
together with $(x,z)$, is $X_0\times\{z\}$. These descriptions
show that isometries preserve both factor foliations.
Every isometry therefore acts separately on $X_0$ and $\R^k$.
Let $E$ be the Euclidean apartment model of $X_0$, and let
$W_0\leq O(E)$ be its reduced finite reflection group.
It is essential, meaning that $E^{W_0}=\{0\}$, or equivalently
that its reflecting normals span $E$. Fix a closed model Weyl
chamber $\mathcal C\subseteq E$, and put
\[
D_0=\{a\in O(E):a\mathcal C=\mathcal C\},\qquad
\Omega_0=W_0\rtimes D_0.
\]
The group $D_0$ is finite: it permutes the finitely many extreme
rays of the essential simplicial cone $\mathcal C$, and an
orthogonal map fixing those rays is the identity. It normalizes
$W_0$ because it permutes the walls of $\mathcal C$ and their
reflections. In particular, $\Omega_0$ is finite.
For $v\in E$ and a unit vector $u\in E$, write
\[
\operatorname{type}(v)=W_0v,\qquad
\operatorname{subtype}(u)=\Omega_0u.
\]
Metric isometries preserve apartments and the reduced Weyl
chambers \cite[Propositions 2.25 and 2.27]{Par00}; on types they
act through $D_0$, and therefore they preserve subtypes.
This is the type/subtype distinction used by Ramos-Cuevas
\cite[Section 2]{RC15}. The residues at individual points need
not be thick.
A building with an entirely thin boundary is a Euclidean
apartment with a chosen Weyl structure. In this case $X_0$
is a point, $E=\{0\}$, and $W_0=\Omega_0=\{I\}$; the
finite-linear-image hypothesis concerns the whole Euclidean
factor. A non-thick building need not be
Euclidean, as $T\times\R$ shows for a branching tree $T$.
Nor must a thin subdivision create a Euclidean factor:
Schwer and Weniger's diagonal subdivision of a product of two
thick trees gives a non-thick $\widetilde B_2$ presentation
\cite[Example 1.1]{SW19}. On the Euclidean factor we allow
arbitrary Euclidean isometries. Finiteness of their linear
image will always be stated as a separate hypothesis.

\begin{remark}[Buildings and polyhedral complexes]\label{rem:buildings-polyhedral}
Discrete affine buildings have piecewise Euclidean Coxeter cell
structures with finitely many cell shapes. General $\R$-buildings
need not have such a structure, or even be locally conical.
For example, attach an infinite half-line to $\R$ at each $1/n$,
$n\geq1$. The resulting complete, geodesically complete $\R$-tree
is a rank-one $\R$-building. Its tangent cone at $0$ is a line,
but every neighbourhood of $0$ contains branch points.
Thus Corollary~\ref{cor:conical-gap} does not apply to all
$\R$-buildings. The general tangent-cone criterion applies when the building
is complete. The proof of Proposition~\ref{prop:building-gap}
instead uses individual ray germs and projections inside
Euclidean apartments, so it also covers incomplete buildings.
\end{remark}

Let $r$ be the apartment dimension of an $\mathbb{R}$-building $X$. Struyve showed \cite[Lemma 4.4]{Str11} that Kleiner and Leeb's result \cite[Theorem 5.1.1]{KL97} that if $X$ is complete it embeds isometrically in a complete $\mathbb{R}$-building $\hat{X}$ of the same type also holds for incomplete $\mathbb{R}$-buildings. Section \ref{sec:dim-ultra} gives $\dim_g \hat{X}= r$. Hence every compact subset $K\subset X$ has covering dimension at most
$r$. A closed ball of positive radius in an apartment is a compact subset of covering
dimension $r$. Therefore
\[\sup\{\dim_{\mathrm{cov}} K : K\subset X \mathrm{compact}\}= r\]
without a completeness assumption on $X$.

Parreau's apartment-projection argument
\cite[Section 4.2]{Par00} already gives the elliptic gap in a
possibly incomplete building; see also Ramos-Cuevas
\cite[Proposition 2.1]{RC15}. We give a self-contained proof
with an explicit, deliberately coarse constant. First we record
the local angular facts, then a quantitative finite-group
lemma, and finally apply them to the fixed set inside an
apartment. Semisimplicity in the original incomplete building
is supplied separately by Le Bars, L\'ecureux and Schillewaert
\cite[Theorem E]{LLS24}.

\begin{lemma}[Local angular facts]\label{lem:building-angles}
In the reduced building $X_0$, two ray germs based at the same
point lie in a common spherical apartment of the residue.
For a unit vector $v\in E$, the angles between germs of subtype
$\Omega_0v$ belong to
\[
 \{\angle(v,wv):w\in\Omega_0\}.
\]
Germs with angle zero have representatives agreeing on an
initial segment. No metric completeness is required.
\end{lemma}
\begin{proof}
Parreau's property (GG) and the residue-building construction
\cite[Corollary 1.11 and Proposition 1.14]{Par00} place the two germs
in a common residue apartment. In a common apartment, the two
directions are represented by $w_1v,w_2v$ with $w_1,w_2\in\Omega_0$,
which gives the angle assertion.
Initial coincidence at angle zero follows from the local
flat-triangle property \cite[Lemma 2.22]{Par00}.
\end{proof}

A fixed subtype has only finitely many possible mutual angles,
but the smallest positive angle can vary with the subtype.
The next lemma perturbs the direction so that its separation
from nonfixed translates dominates the perturbation itself.

\begin{lemma}[Finite orthogonal orbits]\label{lem:finite-angular}
Let $\Omega\leq O(E)$ be finite, put $N=|\Omega|$,
$a=\max\{16,2N\}$, and $\varepsilon=(8a^N)^{-1}$.
For every unit vector $v\in E$, there is a unit vector $u$
with $\angle(u,v)<\pi/2$ such that
\[
 wu\ne u\quad\Longrightarrow\quad
 \angle(u,wu)>2\angle(u,v)+\varepsilon
 \qquad(w\in\Omega).
\]
\end{lemma}
\begin{proof}
Put $\varepsilon_i=a^i\varepsilon$ for $0\leq i\leq N$.
Among the $N$ disjoint intervals
$(\varepsilon_i,\varepsilon_{i+1}]$, one contains none of the
$N-1$ numbers $\|wv-v\|$ with $w\ne I$; fix its index $i$.
Let $S=\{w:\|wv-v\|\leq\varepsilon_i\}$ and $H=\langle S\rangle$.
The set $S$ is symmetric. Every element of $H$ is a word of length
at most $N-1$ in $S$, so it moves $v$ by at most
$(N-1)\varepsilon_i<\varepsilon_{i+1}$.
The empty interval forces that displacement to be at most
$\varepsilon_i$, and therefore $H=S$.

The vector $u_0=|H|^{-1}\sum_{h\in H}hv$ is $H$-fixed and
$\|u_0-v\|\leq\varepsilon_i\leq1/(8a)$.
Hence $u=u_0/\|u_0\|$ is defined,
$\|u-v\|\leq2\varepsilon_i$, and
$\angle(u,v)\leq\pi\varepsilon_i<\pi/2$.
If $wu\ne u$, then $w\notin H$ and
\[
 \|wu-u\|\geq\|wv-v\|-2\|u-v\|
 >(a-4)\varepsilon_i
 >(2\pi+1)\varepsilon_i
 \geq2\angle(u,v)+\varepsilon.
\]
The angle between unit vectors is at least their chord length.
\end{proof}

We now combine the two lemmas with projection onto the
intersection of an apartment and an elliptic fixed set.
Only that intersection needs to be complete. The bound \eqref{eq:building-constant} is deliberately coarse;
no optimality is asserted. For trees the exact constant remains $1/2$.

\begin{proposition}[Building gap without metric completeness]
\label{prop:building-gap}
Let $X=X_0\times\R^k$ be a finite-rank $\R$-building, where
$X_0$ has no Euclidean metric factor. Let a group $G$ act
isometrically on $X$, and suppose that the induced
linear image on $\R^k$ is finite. Then every element of $G$ acts
semisimply, and the action satisfies $\UEG(C)$ for some $C>0$.
When $k=0$, the full isometry group satisfies $\UEG(C_0)$ with
the explicit constant
\begin{equation}\label{eq:building-constant}
 C_0=\frac1{2\sin(\varepsilon/2)},\qquad
 \varepsilon=\frac1{8a^N},\quad
 N=|\Omega_0|,\quad a=\max\{16,2N\}.
\end{equation}
Thus the constant depends only on the canonical reduced spherical
Coxeter complex. Metric completeness is not required.
\end{proposition}

\begin{proof}
We first prove the reduced-factor bound with the constants in
\eqref{eq:building-constant}; the case of a point is immediate.
Let $g$ be elliptic and $x\notin\Fix(g)$.
Choose $q\in\Fix(g)$ and an apartment $A$ containing $q,x$.
The set $F=A\cap\Fix(g)$ is nonempty, closed and convex in the
Euclidean apartment. It is therefore complete, and has a
nearest point $p$ to $x$.
Put $R=d(p,x)>0$ and let $v\in E$ be the unit direction of
$[p,x]$ in a chart of $A$.

Suppose that $\angle_p(x,gx)\leq\varepsilon$.
Apply Lemma~\ref{lem:finite-angular} to $v$ and let $\mu$
be the germ at $p$ of the ray in $A$ with the resulting
direction $u$. The triangle inequality for angles gives
\[
 \angle_p(\mu,g\mu)\leq
 2\angle(u,v)+\angle_p(x,gx)
 \leq2\angle(u,v)+\varepsilon.
\]
Both germs have subtype $\Omega_0u$.
Lemmas~\ref{lem:building-angles} and \ref{lem:finite-angular}
force $\angle_p(\mu,g\mu)=0$.
Thus $g$ fixes an initial segment of this ray.
Since $\langle u,v\rangle>0$, for sufficiently small $t>0$
that segment contains a point $r(t)\in F$ with
\[
 d(x,r(t))^2=R^2+t^2-2tR\langle u,v\rangle<R^2,
\]
contradicting the definition of $p$.
Consequently $\angle_p(x,gx)>\varepsilon$, and \CAT(0) comparison
in the triangle $(p,x,gx)$ yields
\[
 d(x,\Fix(g))\leq R
 \leq\frac{d(x,gx)}{2\sin(\varepsilon/2)}.
\]

Let $H\leq O(k)$ be the finite linear image on $\R^k$.
Choose $C_{\mathrm e}>0$ bounding the reciprocals of the nonzero
singular values of $I-A$ for all $A\in H$. Translating a fixed
point to the origin shows that every elliptic Euclidean isometry
$h$ with linear part in $H$ satisfies
$d(z,\Fix(h))\leq C_{\mathrm e}d(z,hz)$; when its linear part is
the identity, ellipticity makes $h$ the identity.
For an elliptic product isometry the fixed set is the product of
the factor fixed sets, so the squared product-distance formula
gives $\UEG(\max\{C_0,C_{\mathrm e}\})$ on $X$.

Finally, isometries of possibly incomplete $\R$-buildings are
semisimple by \cite[Theorem E]{LLS24}, extending the complete
case of \cite[Th\'eor\`eme 4.1 and Corollaire 4.2]{Par00}.
Apply this on $X_0$. Every Euclidean isometry is semisimple,
and the minimum of the squared displacement of a product is
attained by choosing a minimum point in each factor.
Thus every element of $G$ is semisimple in its action on $X$.
\end{proof}

The passage from a complete enlargement back to the original
building is the same as in the paper of the author, Struyve and Thomas
\cite[Lemma 3.6]{SST25}; we recall it to indicate where finite
generation and the uniform gap enter.
\begin{proof}[Proof of Theorem~\ref{thm:buildings}]
By Proposition~\ref{prop:building-gap} the action on \(X\) satisfies \(\UEG(C)\) for some
\(C>0\). Fix \(o\in X\) and a nonprincipal ultrafilter \(\omega\).
The constant-scale ultrapower construction in
\cite[Lemma 4.4]{Str11} gives a complete building
\[
 \widehat X=\lim_\omega(X,d,o)
\]
of the same finite rank, containing \(X\) equivariantly by constant
sequences. Completeness here is supplied by the ultralimit construction
\cite[Lemma 2.4.2]{KL97}; completeness of $X$ is not needed in
the building construction of \cite[Theorem 5.1.1]{KL97}, as
explained by Struyve \cite[Lemma 4.4]{Str11}.
The extension is explicitly $g[(x_i)]=[(gx_i)]$; it is well
defined since $d(gx_i,o)\leq d(x_i,o)+d(go,o)$.
For each $g$, the basepoint displacement is the constant
$d(o,go)$. Since the original action is locally elliptic,
Lemma~\ref{lem:permanence}(2) makes the induced action on \(\widehat X\) locally
elliptic with the same gap. Since \(\widehat X\) is complete and
has finite geometric dimension, Theorem~\ref{thm:main} gives
\(y\in\Fix_{\widehat X}(G)\). For \(x\in X\),
\[
 d(x,gx)\leq d(x,y)+d(y,gx)=2d(x,y)
 \qquad(g\in G).
\]
Thus \(G\) has a bounded orbit in \(X\). Struyve's fixed-point theorem
for finitely generated bounded groups of isometries of an
\(\R\)-building \cite[Main Result 1]{Str11} gives a fixed point in
\(X\).
\end{proof}

\begin{lemma}[The Euclidean obstruction]\label{lem:euclidean-obstruction}
Let a finitely generated group $G$ act isometrically and locally
elliptically on a finite-rank $\R$-building $X=X_0\times\R^k$,
where $\R^k$ is its Euclidean metric factor. Then the
action on \(X_0\) has a global fixed point, and
\[
 \Fix_X(G)\ne\varnothing
 \quad\Longleftrightarrow\quad
 \Fix_{\R^k}(G)\ne\varnothing.
\]
\end{lemma}
\begin{proof}
Local ellipticity on the product implies local ellipticity on
each factor. Apply Theorem~\ref{thm:buildings} in its reduced case
to \(X_0\). The fixed set on a factor-preserving product is the
product of the group fixed sets, giving the equivalence.
\end{proof}

\begin{remark}[Low-dimensional Euclidean factors]
If $k\leq3$, the action in Lemma~\ref{lem:euclidean-obstruction}
has a global fixed point. The Euclidean assertion used here
requires neither finite generation nor a uniform gap.
For $k=1$, two different reflections generate a nontrivial
translation, so a locally elliptic group of line isometries
has a common fixed point. For $k=2,3$, this is the result of
Breuillard and Fujiwara \cite[Proposition 3.2, proved as
Proposition 9.2]{BF21}; $k=0$ is immediate.
Bass's example in $\R^4$, extended trivially in additional
coordinates, shows that this threshold cannot be increased.
\end{remark}

\begin{remark}[The gap on the Euclidean factor]
For a locally elliptic action on $X=X_0\times\R^k$, the gap has
a stronger Euclidean requirement than the common-fixed-point
conclusion: $\UEG$ holds if and only if the Euclidean affine image
is finite, equivalently if and only if its linear image is finite.
Indeed, given $g$, choose $p_g\in\Fix_{X_0}(g)$ and apply the
product gap at $(p_g,z)$ to obtain the same gap on $\R^k$.
Lemma~\ref{lem:euclidean} makes the Euclidean affine image finite
(the case $k=0$ is immediate). Conversely, finite linear image
gives the gap by Proposition~\ref{prop:building-gap}.
The two finiteness conditions are equivalent because the kernel
of the Euclidean linear-part map consists of elliptic translations,
hence is trivial. This observation does not require finite
generation.
\end{remark}

\begin{lemma}[Finite linear action on a flat]\label{lem:weyl}
Let $X$ be a finite-rank $\R$-building, and let $G$ act
isometrically on $X$. Suppose that the induced linear image
on the Euclidean metric factor
is finite. If a flat $F\subseteq X$ is $G$-invariant, then the
linear image of the action on $F$ is finite. If $\dim F>0$, the
action on $\bd X$ has a finite orbit. Metric completeness of $X$
is not required.
\end{lemma}

\begin{proof}
Write $X=X_0\times\R^k$, and let $E$ be the Euclidean apartment
model of the reduced factor $X_0$. Let $H\leq O(k)$ be the
finite linear image of the action on $\R^k$. With $\Omega_0$
as defined at the beginning of this section, put
\[
M=E\oplus\R^k,\qquad \Omega=\Omega_0\times H\leq O(M).
\]
Thus $\Omega$ is finite.
By Parreau \cite[Proposition 2.25]{Par00}, the flat $F$ is
contained in an apartment $A$. Fix a chart $\phi:A\to M$,
a point $q\in F$, and a linear subspace $V\leq M$ such that
$\phi(F)=\phi(q)+V$. If $V=\{0\}$, the conclusion about the
linear image is immediate, so assume $V\neq\{0\}$.
For each $g\in G$, the restriction of $g$ to $F$ is an affine
Euclidean isometry. Consequently there are unique
$A_g\in O(V)$ and $b_g\in V$ such that
\[
\phi\bigl(g\phi^{-1}(\phi(q)+v)\bigr)
=\phi(q)+b_g+A_gv\qquad(v\in V).
\]
Translations act trivially on the visual boundary of $F$.
It therefore suffices to prove that the set of maps $A_g$
is finite.
On the reduced factor, $g$ preserves segment types up to a
chamber symmetry in $D_0$; on $\R^k$, segment vectors are
transformed by an element of $H$. Applying these two facts
to segments in $F$ gives
\[
A_gv\in\Omega v\qquad(v\in V).
\]
Let $\iota:V\hookrightarrow M$ be the inclusion. Then
\[
V=\bigcup_{w\in\Omega}
\ker(\iota\circ A_g-w\circ\iota).
\]
Each difference in this formula is a linear map $V\to M$.
A real vector space cannot be covered by finitely many proper
linear subspaces. Hence, for some $w\in\Omega$,
$\iota\circ A_g=w\circ\iota$. In particular, $wV=V$ and
$A_g=w|_V$. There are only finitely many such restrictions,
so the linear image of the action on $F$ is finite.
Every unit direction in $V$ therefore has a finite $G$-orbit
in $\partial_\infty F\subseteq\partial_\infty X$.
\end{proof}

\begin{corollary}[Axial random elements]
\label{cor:axial}
Let $X$ be a finite-rank $\R$-building without Euclidean metric
factors, and let a countable group $G$ act isometrically and
non-elementarily on $X$. Let $\mu$ be a symmetric
generating probability measure on $G$ with finite second
displacement moment, and let $Z_n=s_1\cdots s_n$, where the
increments are independent with law $\mu$.
Then the drift is positive and
\[
 \mathbb P(Z_n\text{ acts axially on }X)\longrightarrow1.
\]
Metric completeness of $X$ is not required.
\end{corollary}

\begin{proof}
By Struyve \cite[Lemma 4.4]{Str11}, the constant-scale
ultrapower embeds $X$ equivariantly in a complete building
$\widehat X$ with the apartment model $(E,W_0)$ of the reduced
building $X$.
For each $g\in G$, its action on segment types in $X$ is given
by an element of the fixed finite chamber-symmetry group $D_0$.
This type relation passes to the ultralimit. Hence the same
finite group $\Omega_0$ controls segment types for the extended
action on $\widehat X$. The proof of
Lemma~\ref{lem:weyl}, with $M=E$ and $\Omega=\Omega_0$, therefore
applies to every invariant flat in $\widehat X$.

Let $\overline X$ be the closure of $X$ in $\widehat X$.
The isometric embedding of $X$ is convex, so $\overline X$
is an invariant complete convex \CAT(0) subspace of
$\widehat X$. It is the metric completion of $X$.
Furthermore,
\[
 \gdim\overline X\leq\gdim\widehat X=\operatorname{rank}(X),
\]
and hence $\overline X$ has finite telescopic dimension.

By Le Bars, L\'ecureux
and the author \cite[Proposition 2.6]{LLS24}, the inclusion induces
a natural $G$-equivariant identification
\[
 \bd\overline X=\bd X.
\]
Non-elementarity on $X$ therefore excludes a fixed point
in $\bd\overline X$. It also excludes a fixed point
$y\in\overline X$, since for $o\in X$ such a point would give
\[
 d(o,go)\leq 2d(o,y)\qquad(g\in G).
\]
If $\overline X$ contained an invariant positive-dimensional
flat $F$, regard it as a flat in $\widehat X$.
Lemma~\ref{lem:weyl} would give a finite $G$-orbit in
$\bd F\subseteq\bd\overline X=\bd X$, again contradicting
non-elementarity. An invariant zero-dimensional flat is a
fixed point and has already been excluded.

Choose the basepoint $o$ in $X$. For $o\in X$ and every $g\in G$, the equivariant inclusion
$X\hookrightarrow\overline X$ preserves $d(o,go)$. In particular,
it preserves the second displacement moment and, for every
sample path and every $n$, the value $d(o,Z_no)$. The drift
is therefore the same in the two spaces.
Theorem~\ref{thm:izeki},
applied to $\overline X$, now implies that this drift is
positive. On $X$, Proposition~\ref{prop:building-gap}
gives a uniform elliptic gap and semisimplicity.
Corollary~\ref{cor:zero-drift} now gives the asserted
probability of axiality.
\end{proof}

\begin{remark}[Quantitative consequences]
The argument does not itself give a uniform word-length bound
for a nonelliptic element or a rate in
Corollary~\ref{cor:axial}. Any suitable scalar displacement
concentration estimate supplies quantitative conclusions
through Corollary~\ref{cor:concentration}. For Bruhat--Tits
buildings, a quantitative escape theorem is already available
from Breuillard and Fujiwara \cite[Theorem 7.16]{BF21}.
\end{remark}

\begin{remark}[Strong regularity]\label{rem:quint}
An axial building isometry is \emph{strongly regular} if the endpoints
of an axis lie in the interiors of opposite chambers in the spherical
building at infinity. Axiality in Corollary~\ref{cor:axial} cannot in
general be replaced by strong regularity. To see this, let $k$ be a
non-Archimedean local field. The Jordan
projection of an element of $\mathrm{Sp}_4(k)$ is its translation
vector, represented in the closed positive Weyl chamber; the
\emph{limit cone} of a subgroup is the closed cone generated by
its Jordan projections. In type $C_2$, the longest Weyl element is
$-I$, so the opposition involution is the identity. A nonzero rational
ray in a wall is therefore invariant under opposition and is not
contained in the central subspace.

Quint's realization theorem \cite[Th\'eor\`eme 5.2]{Qui02} applies to a cone in the positive chamber that is
“un cône convexe, fermé, à support rationnel stable par $\iota$
et non inclus dans $E_C$”, where $E_C$ is the central subspace.
There is no nonempty-interior hypothesis.
The rational wall ray above satisfies these conditions, so the
theorem gives a discrete Zariski-dense subgroup of $\mathrm{Sp}_4(k)$ with this
limit cone. All its Jordan projections lie in the wall, so it has no
strongly regular elements. Its action on the proper Bruhat--Tits
building is non-elementary. Indeed, a bounded orbit would give a fixed point,
putting the discrete subgroup in a compact stabilizer and making it
finite. A finite boundary orbit would put a finite-index subgroup
in a proper parabolic subgroup, contrary to Zariski density in the
connected algebraic group. See also \cite[Remark 1.3]{LLS24}.
\end{remark}

\section{Finite-dimensional cube complexes}\label{sec:cubes}

Here $X$ is a \CAT(0) cube complex with unit cubes, $X^{(0)}$
is its vertex set, and $d_2$ is its intrinsic piecewise Euclidean
metric. A cubical automorphism restricts on each cube to a Euclidean
cube isometry, hence preserves lengths of piecewise cubical
paths. Applying this to the automorphism and its inverse shows
that it is an isometry of the intrinsic metric $d_2$; the same
argument applies to the intrinsic piecewise $\ell^1$ metric.
The edge-path distance $d_1(u,v)$ between vertices counts the
hyperplanes separating them. We also write $d_1$ for the
intrinsic piecewise $\ell^1$ metric, which restricts to that
distance on vertices; see Haglund \cite[Section 2.2]{Hag23} and
Leder and Varghese \cite[Section 2]{LV20}.
A hyperplane has two vertex halfspaces; two hyperplanes
\emph{cross} if all four intersections of their halfspaces
are nonempty. A pairwise crossing family is dual to a cube.

The full subcomplex determined by a nonempty intersection of
vertex halfspaces has a \emph{gate}: for each vertex $y$ there
is a unique vertex $z$ in the subcomplex such that
\[
 d_1(y,u)=d_1(y,z)+d_1(z,u)
 \quad\text{for every vertex $u$ of the subcomplex}.
\]
In particular, $z$ is the nearest vertex to $y$, see Fern\'os--Futer--Hagen
\cite[Lemma 3.18 and Corollary 3.20]{FFH24}.
If a cubical automorphism $g$ preserves such a subcomplex $Q$,
then $g\pi_Q(y)=\pi_Q(gy)$: applying $g$ to the defining gate
identity shows that $g\pi_Q(y)$ is the gate of $gy$, and gates
are unique.

A connected graph is \emph{median} if any three vertices have
a unique vertex lying on shortest paths between each pair.
Its cubical dimension is the supremum of the dimensions of its
hypercubes. By Chepoi's characterization \cite[Theorem 6.1]{Che00},
median graphs are precisely the vertex graphs of $\CAT(0)$
cube complexes, obtained by filling their hypercubes.

\begin{theorem}[Median-graph elliptic gap]\label{thm:median-gap}
Let $M$ be a median graph of cubical dimension at most $d\geq1$,
with unit edge-path metric $d_1$. If $g\in\operatorname{Aut}(M)$ fixes a vertex, then for every vertex
$x$ there is a fixed vertex $y$ such that
\[
 d_1(x,y)\leq\frac d2\,d_1(x,gx).
\]
\end{theorem}
\begin{proof}
Regard $M$ as the vertex graph of its associated cube complex $X$;
the automorphism $g$ extends cubically.
Choose a fixed vertex \(y\) minimizing the integer \(d_1(x,y)\).
Let \(\mathcal H\) be the finite set of hyperplanes separating \(x,y\).
For every $H\in\mathcal H$, denote by $h_H$ the vertex
halfspace bounded by $H$ and containing $y$; its complementary
halfspace contains $x$. We use these choices to orient the
hyperplanes in $\mathcal H$.

Distinct translates under \(\langle g\rangle\) of a halfspace
containing \(y\) cannot be strictly nested.
Otherwise, after reversing the power if necessary, there is a strict
descending chain
\(h\supsetneq g^kh\supsetneq g^{2k}h\supsetneq\cdots\),
all containing \(y\).
Any vertex outside \(h\) would then be separated from \(y\) by
infinitely many hyperplanes, which is impossible.

It follows that the members of any one hyperplane orbit lying in
\(\mathcal H\) are pairwise crossing: both the \(y\)-sides and the
\(x\)-sides intersect, and neither oriented halfspace contains the other.
Thus each orbit contributes at most \(d\) elements to \(\mathcal H\).

No entire hyperplane orbit can be contained in \(\mathcal H\).
If one were, intersect its \(x\)-halfspaces.
The resulting combinatorially convex subcomplex $Q$ contains \(x\), excludes \(y\),
and is \(g\)-invariant.  The gate \(z\) of \(y\) in \(Q\) is fixed by
\(g\), and
\(d_1(x,z)=d_1(x,y)-d_1(y,z)<d_1(x,y)\), a contradiction.

Hence every orbit meeting \(\mathcal H\) contributes at least one element
to \(\mathcal H\setminus g\mathcal H\).
Because \(gy=y\), the set \(g\mathcal H\) consists of the hyperplanes
separating $gx,y$. Since $|\mathcal H|=|g\mathcal H|$, one obtains
\[
 d_1(x,y)=|\mathcal H|
 \leq d\,|\mathcal H\setminus g\mathcal H|
 =\frac d2\,|\mathcal H\mathbin{\triangle}g\mathcal H|
 =\frac d2\,d_1(x,gx).
\]
\end{proof}

\begin{lemma}[Sparse permutation estimate]
\label{lem:sparse-permutation}
Let $V$ be a set, let $\sigma:V\to V$ be a bijection, and let
$P$ be its permutation operator on the real Hilbert space
$\ell^2(V;\mathbb R)$, defined by
\[
(Pa)(v)=a(\sigma^{-1}v)\qquad(v\in V).
\]
Let $D\geq1$ be an integer and let $a\in\ell^2(V;\mathbb R)$
satisfy $a(v)\geq0$ for every $v\in V$ and
$|\operatorname{supp}(a)|\leq D$. There exists
$b\in\ell^2(V;\mathbb R)$ such that $b(v)\geq0$ for every
$v\in V$, $Pb=b$, $\operatorname{supp}(b)\subseteq
\operatorname{supp}(a)$, and
\[
\|a-b\|_2\leq\beta_D\|a-Pa\|_2,
\qquad
\beta_D=\frac1{2\sin(\pi/(2(D+1)))}.
\]
\end{lemma}
\begin{proof}
On every finite $P$-orbit contained in $\operatorname{supp}(a)$,
replace the coordinates by their mean; set all other coordinates
of $b$ to zero. This defines a fixed vector supported in
$\operatorname{supp}(a)$.
On a full orbit of length $m>1$, the cyclic
inequality of Fan--Taussky--Todd \cite[Theorem 10, p.~80]{FTT55} gives
\[
 \|a-b\|\leq\frac{\|a-Pa\|}{2\sin(\pi/m)}
 \leq\beta_D\|a-Pa\|,
\]
since $m\leq D$.
On any other orbit, the support splits into consecutive blocks
of lengths $r\leq D$, with zero boundary values.
For such a block $(a_1,\ldots,a_r)$, their Dirichlet inequality
\cite[Theorem 9, p.~79]{FTT55} gives
\[
 \sum_{j=0}^{r}(a_{j+1}-a_j)^2
 \geq4\sin^2\!\left(\frac{\pi}{2(r+1)}\right)
       \sum_{j=1}^{r}a_j^2,\qquad a_0=a_{r+1}=0.
\]
Summing the cycle and block estimates proves the assertion.
\end{proof}

\begin{proposition}[Cubical elliptic gap]\label{prop:cube-intro}
Let $X$ be a \CAT(0) cube complex of dimension at most $d\geq1$.
For its unit piecewise Euclidean metric, the action of its cubical
automorphism group satisfies $\UEG(C_d)$, where
\begin{equation}\label{eq:cubical-constant}
 C_1=\frac12,\qquad
 C_d=\frac1{2\sin(\pi/(2(d+1)))}\quad(d\geq2).
\end{equation}
For the intrinsic piecewise $\ell^1$ metric the action satisfies
$\UEG(d/2)$.
\end{proposition}
\begin{proof}
Let $p$ be a fixed point of the elliptic cubical automorphism
$g$, and let $Q$ be the unique cube whose relative interior
contains $p$. Since $g$ is cubical and fixes $p$, it preserves
$Q$ and hence fixes its centre. In the subdivision obtained
by bisecting every edge, that centre is a vertex. Thus $g$
fixes a vertex in this and in every further dyadic subdivision.
Apply Theorem~\ref{thm:median-gap} in every further
dyadic subdivision, with its original edge lengths.
The dyadic vertices are dense and distance to the fixed set is
$1$-Lipschitz, so the $\ell^1$ estimate passes to all points.

For the Euclidean estimate, $X$ is complete by
\cite[I.7, Theorem 7.50]{BH99}. The case $d=1$ is the tree
bound. Suppose $d\geq2$, let $g$ fix $p$, and let $Q$ be the
smallest cube containing $p$, of dimension $k$.
The tangent cone has the form
\[
 T_pX=\R^k\times\mathcal C,
\]
where $\mathcal C$ is the union of nonnegative coordinate
orthants indexed by simplices in the link of $Q$.
Each vector in $\mathcal C$ has at most $d-k$ nonzero
coordinates. The action $D_pg$ is a signed permutation on
$\R^k$ and a permutation of the coordinates of $\mathcal C$.
A neighbourhood of $p$ is isometric to a neighbourhood of
the cone vertex \cite[I.7, Theorem 7.16]{BH99};
in particular, every $D_pg$-fixed vector is tangent to an
actual fixed segment.

Embed these coordinate orthants in $\ell^2(V)$.
The coordinate map $\mathcal C\to\ell^2(V)$ is $1$-Lipschitz:
it preserves lengths within each orthant, while the cone metric
is intrinsic. Distances between vectors in one orthant agree
with their Euclidean coordinate distances.
For $a\in\mathcal C$, Lemma~\ref{lem:sparse-permutation}
gives a fixed vector $b$ in the same orthant with
$\|a-b\|\leq\beta_d\|a-Pa\|$; the zero vector needs no estimate.

On $\R^k$, let $A$ be the signed permutation and $z_0$ the
orthogonal projection of $z$ onto $\ker(I-A)$.
A signed cycle of length $j\leq k$ has eigenvalues satisfying
$\lambda^j=1$ or $\lambda^j=-1$.
Consequently its nonidentity eigenvalues have angular distance
at least $\pi/k$ from $1$, and, when $k>0$,
\[
 \|z-z_0\|\leq
 \frac{\|z-Az\|}{2\sin(\pi/(2k))}
 \leq\beta_d\|z-Az\|.
\]
The vector $(z_0,b)$ is tangent to $\Fix_X(g)$.
The product distance formula and the coordinate estimates give
\[
 d_p\bigl((z,a),T_p\Fix_X(g)\bigr)
 \leq\beta_d\,d_p\bigl((z,a),D_pg(z,a)\bigr).
\]
Theorem~\ref{thm:tangent-gap} proves the Euclidean bound.
Finally, concavity of sine on $[0,\pi/6]$ gives
$\sin(\pi/(2(d+1)))\geq3/(2(d+1))\geq1/d$ for $d\geq2$,
so $C_d\leq d/2$.
\end{proof}

\begin{remark}[Optimality of the constants]\label{rem:cubical-sharpness}
The coefficient $d/2$ in Theorem~\ref{thm:median-gap} and in
the $\ell^1$ estimate is optimal.
On the standard cubulation of $\R^d$, the cubical isometry
$A(x_1,\ldots,x_d)=(-x_d,x_1,\ldots,x_{d-1})$ fixes only $0$.
For $u=(1,\ldots,1)$, one has $d_1(u,0)=d$ and $d_1(u,Au)=2$.
For the Euclidean metric this example has optimal constant
$(2\sin(\pi/(2d)))^{-1}$, since the eigenvalues of $A$
are the roots of $z^d=-1$.
Thus, if $C_d^*$ denotes the optimal universal Euclidean
constant, then $C_1^*=1/2$ and, for $d\geq2$,
\[
 \frac1{2\sin(\pi/(2d))}\leq C_d^*
 \leq\frac1{2\sin(\pi/(2(d+1)))}.
\]
In particular, $C_d^*\sim d/\pi$.
The exact optimum in each higher dimension is not asserted.
\end{remark}

\begin{corollary}[Cubical local-to-global fixed points]\label{cor:cubical-fixed}
A finitely generated group acting locally elliptically by cubical
automorphisms on a finite-dimensional \CAT(0) cube complex fixes
a point.
\end{corollary}
\begin{proof}
The unit piecewise Euclidean metric is complete and has finite
geometric dimension. Apply Proposition~\ref{prop:cube-intro}
and Theorem~\ref{thm:main}.
\end{proof}
The qualitative fixed-point conclusion is due to Sageev
\cite[proof of Theorem 5.1]{Sag95} and was made explicit by
Leder--Varghese \cite[Theorem A]{LV20}. The contribution here
is a quantitative displacement bound: the constant $d/2$ is
sharp for the median-graph and piecewise $\ell^1$ metrics,
and the universal Euclidean constant has the sharp asymptotic
growth $d/\pi$. Genevois--Lonjou--Urech
\cite[Theorem 1.3]{GLU24} prove the qualitative conclusion for
$\CAT(0)$ cube complexes without infinite cubes, without
requiring a uniform bound on their dimensions.

\section{Further consequences of the elliptic gap}\label{sec:examples}

\subsection{Finite orders and bounded stabilizers}

Finite-order isometries satisfy a sharp displacement estimate
on every complete $\CAT(0)$ space. The proof combines the
variance inequality for a finite orbit with a nonlinear
Wirtinger inequality.

\begin{proposition}\label{prop:order}
Let $X$ be a complete \CAT(0) space, and let $g\in\operatorname{Isom}(X)$
satisfy $g^m=1$ for an integer $m\geq2$. Then
\[
 d\bigl(x,\operatorname{Fix}(g)\bigr)\leq\frac{1}{2\sin(\pi/m)}\,d(x,gx)
 \qquad(x\in X).
\]
For an isometry of order $m$, this constant is optimal over complete
\CAT(0) spaces.
\end{proposition}

\begin{proof}
Fix $x$, put $L=d(x,gx)$, and consider
\[
 F(y)=\frac1m\sum_{i=0}^{m-1}d(y,g^ix)^2 .
\]
Completeness and strong \CAT(0) convexity give a unique minimizer $p$.
Invariance of $F$ under $g$ gives $gp=p$. Put $R=d(p,x)$; then
$F(p)=R^2$.

For $p_t$ at parameter $t$ on $[p,y]$, strong convexity and minimality give
\[
 F(p)\leq F(p_t)\leq(1-t)F(p)+tF(y)-t(1-t)\,d(p,y)^2 .
\]
After division by $t$ and passage to $t\downarrow0$, this gives
$F(y)\geq F(p)+d(p,y)^2$. In particular, for $y=x$,
\begin{equation}\label{eq:R-bound}
 2R^2\leq\frac1m\sum_{j=0}^{m-1}d(x,g^jx)^2 .
\end{equation}

For $m\geq4$, the nonlinear Wirtinger inequality for \CAT(0) spaces,
applied to the cyclic list $g^ix$, states that
\[
 \sin^2(\pi/m)\sum_{i=0}^{m-1}d(g^ix,g^{i+j}x)^2
 \leq\sin^2(j\pi/m)\sum_{i=0}^{m-1}d(g^ix,g^{i+1}x)^2
\]
for $2\leq j\leq m-2$. 

We use Gromov's nonlinear Wirtinger inequality for $\CAT(0)$
spaces, in the formulation recorded by Toyoda
\cite[Section 1.2, Definition 1.7 and equation (1.2)]{Toy23}.
Cyclic symmetry gives
\[
 d(x,g^jx)^2\leq\frac{\sin^2(j\pi/m)}{\sin^2(\pi/m)}\,L^2 .
\]
The same inequality for $j=0,1,m-1$ is immediate. Summing it and using
$\sum_{j=0}^{m-1}\sin^2(j\pi/m)=m/2$ in \eqref{eq:R-bound} gives
\[
 2R^2\leq\frac{L^2}{2\sin^2(\pi/m)},
\]
that is, $R\leq L/(2\sin(\pi/m))$. Since $p\in\operatorname{Fix}(g)$, we
have $d(x,\operatorname{Fix}(g))\leq R$, which proves the estimate. For
$m=2,3$, all the needed summands are among the immediate cases, so the
same argument works without invoking a nontrivial Wirtinger inequality.

\emph{Optimality.} Let $X=\mathbb R^2$ and let $g$ be the rotation about the
origin through the angle $2\pi/m$. Then $g$ has order exactly $m$ and
$\operatorname{Fix}(g)=\{0\}$. For $x\in\mathbb R^2$ with $|x|=r>0$, the
points $x$ and $gx$ span a chord of the circle of radius $r$ subtending the
angle $2\pi/m$, so
\[
 d(x,gx)=2r\sin(\pi/m),
 \qquad
 d\bigl(x,\operatorname{Fix}(g)\bigr)=r=\frac{d(x,gx)}{2\sin(\pi/m)} .
\]
Thus equality holds in the estimate, and the constant cannot be improved.
\end{proof}

Orders in Proposition~\ref{prop:order} are orders of isometries in the image of the action. For each fixed $N$, the following consequence applies the estimate only to those values of $Z_n$ whose induced isometry has order at most $N$, using the block interchange of Hou--Izeki--Ji--Karlsson--Wu
\cite[Lemma 3.2]{HIJKW26}; no uniform elliptic gap is assumed.

\begin{corollary}[Escape from bounded orders]\label{cor:orders}
Let a countable group $G$ act isometrically on a \CAT(0) space $(X,d)$,
let $o\in X$, and let $\mu$ be a probability measure on $G$ with
$\sum_g\mu(g)d(o,go)<\infty$.
For independent increments $s_i$ of law $\mu$, let
$Z_n=s_1\cdots s_n$. If the drift is positive, then, for every
integer $N\geq1$,
\[
 \mathbb P\bigl(\text{the isometry induced by }Z_n
                       \text{ has order at most }N\bigr)
 \longrightarrow0.
\]
\end{corollary}
\begin{proof}
If $N=1$, the event in question implies $d(o,Z_no)=0$,
whose probability tends to zero by positive drift.
Assume $N\geq2$ and extend the action to the metric completion
$\overline X$, which is complete and $\CAT(0)$.
Displacements at $o$, the drift, and the exact orders of the
induced isometries are unchanged. On the event in question,
Proposition~\ref{prop:order} in $\overline X$ and
Lemma~\ref{lem:ratio-gap} apply with
\[
C=\frac1{2\sin(\pi/N)}.
\]
They give a uniform upper bound strictly below $2$ for
$d(o,Z_n^2o)/d(o,Z_no)$ when the denominator is nonzero.
Proposition~\ref{prop:ratio} makes the probability of this
event tend to zero; positive drift also gives
$\mathbb P(d(o,Z_no)=0)\to0$.
\end{proof}

\begin{corollary}[Uniformly bounded point stabilizers]\label{cor:stabilizers}
Let a group $\Gamma$ act isometrically on a complete \CAT(0) space
$X$.  Suppose that there is an integer
\(N\geq1\) such that
\[
 |\Gamma_x|\leq N\quad(x\in X),\qquad
 \Gamma_x=\{\gamma\in\Gamma:\gamma x=x\}.
\]
If $N\geq2$, the action satisfies
$\UEG((2\sin(\pi/N))^{-1})$.
If $N=1$, every positive gap constant is valid.
If $X$ has finite geometric dimension, every torsion subgroup
of $\Gamma$ has at most $N$ elements.
\end{corollary}
\begin{proof}
An elliptic element lies in a point stabilizer, so its order is at most
\(N\).  Proposition~\ref{prop:order} gives the gap.
Suppose in addition that $X$ has finite geometric dimension.
Every finitely generated torsion subgroup acts locally elliptically
and fixes a point by Theorem~\ref{thm:main}; it therefore has at most
$N$ elements. If a torsion subgroup contained $N+1$ distinct
elements, the subgroup they generate would contradict this bound.
\end{proof}

The torsion-subgroup conclusion also follows from
Hou--Izeki--Ji--Karlsson--Wu
\cite[Theorem 1.2(2)]{HIJKW26}: every element of a torsion
subgroup has order at most $N$, hence its exponent divides
$\operatorname{lcm}(1,\ldots,N)$.
The displacement estimate above is sharp as a universal
bounded-order estimate, and Corollary~\ref{cor:orders} applies
it probabilistically without any dimension assumption.
For faithful actions on Hadamard manifolds, uniformly bounded
finite point stabilizers are also necessary for a uniform gap
(Proposition~\ref{prop:smooth}).

\subsection{Geometric actions}
An action of a discrete group \(\Gamma\) on a proper \CAT(0) space
\(X\) is \emph{geometric} if it is proper and cocompact.  Here
properness means that
\(\{\gamma\in\Gamma:\gamma K\cap K\neq\varnothing\}\) is finite for
every compact subset \(K\subseteq X\), and cocompactness means that
\(\Gamma K=X\) for some compact \(K\).
Such an action has uniformly bounded finite point stabilizers: move
each point into \(K\), and observe that its stabilizer is contained in
the finite set just displayed.  Proposition~\ref{prop:order} therefore
gives a uniform elliptic gap.

\begin{corollary}[Random axial elements for geometric actions]
\label{cor:geometric}
Let a nonamenable discrete group $\Gamma$ act geometrically
on a proper \CAT(0) space $X$. Let $\mu$ be a symmetric
generating probability measure on $\Gamma$ such that
\[
\sum_{\gamma\in\Gamma}\mu(\gamma)d(o,\gamma o)<\infty
\]
for some $o\in X$. For independent increments $s_i$ of law
$\mu$, put $Z_n=s_1\cdots s_n$. Then the drift is positive and
\[
\mathbb P(Z_n\text{ acts axially on }X)\longrightarrow1.
\]
\end{corollary}
\begin{proof}
Put $L(\gamma)=d(o,\gamma o)$. The \v{S}varc--Milnor lemma
\cite[I.8, Proposition 8.19]{BH99} gives a finite generating
set and a quasi-isometric orbit map. Consequently there are
$A,a>0$ such that
\[
B_L(R)=\{\gamma\in\Gamma:L(\gamma)\leq R\},\qquad
|B_L(R)|\leq Ae^{aR}\quad(R\geq0).
\]
This counts elements of $\Gamma$, including the finite
multiplicity caused by a point stabilizer.
Let $P_\mu$ be the convolution operator on $\ell^2(\Gamma)$.
Kesten's criterion \cite{Kes59a,Kes59b} gives
$\rho=\|P_\mu\|<1$. The criterion applies to infinite-support
measures as well; see also \cite[Section 1]{DS24}.
For completeness, choose a finite symmetric generating
subset $F$ of $\operatorname{supp}\mu$, put
$\alpha=\mu(F)>0$ and $\nu=\alpha^{-1}\mu|_F$.
If $\alpha<1$, write
$\mu=\alpha\nu+(1-\alpha)\eta$. The finite-support criterion
and $\|P_\eta\|\leq1$ give
\[
\|P_\mu\|\leq\alpha\|P_\nu\|+1-\alpha<1.
\]
If $\alpha=1$, apply the finite-support criterion directly.
Cauchy--Schwarz now gives
\[
\mathbb P(L(Z_n)\leq R)
\leq |B_L(R)|^{1/2}\|\mu^{*n}\|_2
\leq A^{1/2}e^{aR/2}\rho^n.
\]
Choose $c>0$ such that $e^{ac/2}\rho<1$.
The probabilities of $L(Z_n)\leq cn$ are summable, so
Borel--Cantelli gives
$\liminf_n L(Z_n)/n\geq c$ almost surely.
The first moment and Kingman's theorem imply that the
drift exists and satisfies $\ell\geq c>0$.
Proper cocompactness gives a uniform bound on finite point
stabilizers, so Proposition~\ref{prop:order} supplies a uniform
elliptic gap. Every acting isometry is semisimple
\cite[II.6, Proposition 6.10(2)]{BH99}.
The conclusion follows from Corollary~\ref{cor:zero-drift}.
\end{proof}

\begin{remark}[Non-elementary geometric actions]
The conclusion applies, in particular, to every non-elementary
geometric action. Indeed, if $\Gamma$ were amenable, the
Adams--Ballmann theorem \cite{AB98} would give a fixed point
in $\partial_\infty X$ or an invariant Euclidean flat
(allowing a point); see also \cite[Theorem 5.1]{Cap09}.
The point and boundary alternatives contradict
non-elementarity. If $F$ is an invariant positive-dimensional
flat and $\Gamma K=X$ with $K$ compact, then
$\Gamma(K\cap F)=F$. Thus the restricted action is proper
and cocompact. Its kernel is finite, and its Euclidean
isometry image is crystallographic. Bieberbach's theorem
\cite{Bie11,Bie12} gives finite linear image, hence a finite
orbit in $\partial_\infty F\subseteq\partial_\infty X$,
again a contradiction.
The nonamenable formulation is strictly more general.
For example, $F_2\times\mathbb Z$ acts geometrically on
$T\times\R$, where $T$ is the Cayley tree of $F_2$.
The action fixes the two endpoints of the $\R$-factor and
is therefore not non-elementary, but the corollary applies.
\end{remark}

An axial isometry is \emph{rank one} if none of its axes bounds a
Euclidean half-plane. A geodesic line is \emph{contracting} if there
is a uniform bound on the diameters of projections onto it of metric
balls disjoint from it; an axial isometry is contracting if it has
such an axis. For non-elementary actions on proper \CAT(0) spaces
containing a rank-one isometry, Le Bars proves that random elements
are contracting with probability tending to one, for every measure
whose support generates the group as a semigroup
\cite[Theorem 1.2]{LB26}. That theorem requires neither symmetry nor
a moment condition. Corollary~\ref{cor:geometric} assumes symmetry and a finite
first displacement moment, but requires no rank-one element
and applies also to geometric nonamenable actions with finite
boundary orbits. Its conclusion is axiality.

\subsection{Products and subgroups}

Convex restrictions preserve the gap by Lemma~\ref{lem:permanence}.
For products, let $G$ act isometrically on $\CAT(0)$ spaces $X_1,\dots,X_r$, the $i$-th action satisfying $\UEG(C_i)$, and give $X=\prod_iX_i$ the metric $d^2=\sum_id_i^2$ and the diagonal action. If $g$ is elliptic, then $\Fix_X(g)=\prod_i\Fix_{X_i}(g)$, and since the infima are independent (no completeness is needed),
\[
 d(x,\Fix_X(g))^2=\sum_i d_i(x_i,\Fix_{X_i}(g))^2\leq(\max_iC_i)^2\,d(x,gx)^2,
\]
so the action on $X$ satisfies $\UEG(\max_iC_i)$.

\begin{proposition}[Permanence under finite-index overgroups]
\label{prop:finite-index}
Let \(G\) act isometrically on a complete \CAT(0) space \(X\), and
let \(H\leq G\) have finite index \(m\geq 2\). If the restricted action of
\(H\) satisfies \(\UEG(C)\), then the action of \(G\) satisfies
\[
 \UEG\bigl(mC+\frac{1}{2\sin(\pi/m)}).
\]
\end{proposition}
\begin{proof}
Let $g\in G$ be elliptic. If $g\in H$, then $d(x,\Fix(g))\leq C\,d(x,gx)$
by assumption, so assume $g\notin H$. Among the $m+1$ cosets
$H,gH,\ldots,g^mH$ two coincide, so $g^k\in H$ for some $2\leq k\leq m$.
The set \(D=\Fix(g^k)\) is nonempty,
complete, convex, and \(g\)-invariant. Let \(p\) be the projection
of \(x\) onto \(D\). Then
\[
 d(x,p)\leq C\,d(x,g^kx)\leq Ck\,d(x,gx).
\]
Projection commutes with \(g\), so \(d(p,gp)\leq d(x,gx)\).
The restriction of \(g\) to \(D\) has order dividing \(k\).
Proposition~\ref{prop:order}, applied in \(D\), gives a fixed point \(q\) with
\[
 d(p,q)\leq \frac{1}{2\sin(\pi/k)}\,d(p,gp).
\]
Since \(q\) is fixed by \(g\) in \(X\), the triangle inequality
gives the result.
\end{proof}

\section{Hadamard manifolds}\label{sec:smooth}

A \emph{Hadamard manifold} $M$ is a complete simply connected Riemannian manifold of nonpositive sectional curvature; it is a proper $\CAT(0)$ space \cite[II.1A and the Cartan--Hadamard theorem II.4.1]{BH99}, and $\exp_p\colon T_pM\to M$ is a diffeomorphism for each $p\in M$. Isometries of $M$ are smooth by the Myers--Steenrod theorem \cite{MSt39}. If $g$ fixes $p$, then $A=D_pg$ is orthogonal and, writing $E=\ker(I-A)$,
\[
 g\exp_p(u)=\exp_p(Au),\qquad \Fix_M(g)=\exp_p(E).
\]
Hence $D_pg=I$ forces $g=1$, so $D_p\colon\Gamma_p\to O(T_pM)$ is injective for every subgroup $\Gamma\leq\Isom(M)$ with $\Gamma_p=\{g\in\Gamma:gp=p\}$.
The subspace $E$ is the tangent space at $p$ of the closed convex set $\Fix_M(g)$, and its orthogonal complement $E^\perp$ is the \emph{normal space}. Since the directions at $p$ towards other fixed points are the nonzero $w\in E$, and $E=-E$, the $\CAT(0)$ projection criterion (angle at least $\pi/2$ with every such $w$) gives
\[
 p=\operatorname{proj}_{\Fix_M(g)}(\exp_p v)\iff v\perp E,
 \qquad\text{and then}\qquad
 d(\exp_p(tv),\Fix_M(g))=t\|v\|\quad(t\geq0).
\]

The hinge comparison inequality
\cite[II.1, Proposition 1.7(5) and Exercise 1.9(2)]{BH99}
gives the first inequality below: the Alexandrov angle at $p$
is the Riemannian angle between $u$ and $w$, and the Euclidean
cosine law gives the lower bound $\|u-w\|$.
The second equality is the first-order expression of the
Riemannian distance in normal coordinates at $p$.

\begin{align}
 d(\exp_p u,\exp_p w)&\geq\|u-w\|,\label{eq:smooth-exp-comparison}\\
 \lim_{t\downarrow0}\frac{d(\exp_p(tu),\exp_p(tw))}{t}&=\|u-w\|.\label{eq:smooth-tangent-metric}
\end{align}

\begin{proposition}[Isotropy criterion]\label{prop:smooth}
Let $M$ be an $n$-dimensional Hadamard manifold, where
$n\geq1$, and let $\Gamma\leq\Isom(M)$. The following
conditions are equivalent:
\begin{enumerate}
\item The action of $\Gamma$ has a uniform elliptic gap.
\item The orders of its elliptic elements are finite and
uniformly bounded.
\item The stabilizers
$\Gamma_p=\{g\in\Gamma:gp=p\}$ are finite, with cardinalities
bounded uniformly for $p\in M$.
\end{enumerate}
For each $C>0$, the action satisfies $\UEG(C)$ if and only
if, for every nonidentity elliptic $g\in\Gamma$, every
$p\in\Fix_M(g)$, and every
$v\in\ker(I-D_pg)^\perp$, one has
\begin{equation}\label{eq:normal-gap}
 \|(I-D_pg)v\|\geq C^{-1}\|v\|.
\end{equation}
Here norms and orthogonal complements are taken in the
Riemannian tangent space. Moreover, $\UEG(C)$ implies
\begin{equation}\label{eq:smooth-bound}
 |\Gamma_p|\leq\mathcal B_n(C)\qquad(p\in M).
\end{equation}
If the identity is the only elliptic element, every $C>0$ works.
If an integer $N\geq2$ bounds the orders in \textup{(2)}
or the cardinalities in \textup{(3)}, then the action satisfies
$\UEG((2\sin(\pi/N))^{-1})$.
\end{proposition}

\begin{proof}
For $p\in\Fix_M(g)$, the tangent cone of the fixed set is
$\ker(I-D_pg)$. Theorem~\ref{thm:tangent-gap} therefore identifies
$\UEG(C)$ with \eqref{eq:normal-gap}.
For each point stabilizer, apply Lemma~\ref{lem:orthogonal}
to the faithful derivative representation
$D_p\Gamma_p\leq O(T_pM)$. This proves
\eqref{eq:smooth-bound} and \textup{(1)}$\Rightarrow$\textup{(3)}.
Every elliptic element belongs to a point stabilizer, so
\textup{(3)} implies \textup{(2)}.

Finally, suppose that \textup{(2)} holds with bound
$N\geq2$. Let $g\ne1$ be elliptic of order $m\leq N$,
choose $p\in\Fix_M(g)$, and put $A=D_pg$.
The orthogonal map $A$ satisfies $A^m=I$.

Since $A$ is orthogonal, its complexification is unitary and
$I-A$ is normal. Since $A^m=I$, its complex eigenvalues
different from $1$ are $m$th roots of unity. Thus the
corresponding nonzero singular values of $I-A$ are
\[
|1-e^{2\pi i k/m}|
=2|\sin(\pi k/m)|
\geq2\sin(\pi/m),
\qquad 1\leq k\leq m-1.
\]

Therefore
\[
 \|(I-A)v\|
 \geq2\sin(\pi/m)\|v\|
 \geq2\sin(\pi/N)\|v\|
 \qquad(v\perp\ker(I-A)).
\]
The normal criterion gives the asserted gap and proves
\textup{(2)}$\Rightarrow$\textup{(1)}.
A cardinality bound $N$ in \textup{(3)} also bounds the
order of each elliptic element, so it gives the same constant.
If $N=1$, the only elliptic element is the identity,
and every positive gap constant is valid.
\end{proof}

\begin{remark}[Small rotations]\label{rem:small-rotations}
The full isometry group of the hyperbolic
plane does not have a uniform elliptic
gap, since rotations about a fixed point
have normal derivatives with arbitrarily
small nonzero rotation angles.
\end{remark}

\begin{remark}[Lattices and their subgroups]\label{rem:lattices}
Let $\Lambda$ be a discrete subgroup of $\Isom(M)$ with a
torsion-free normal subgroup $\Lambda_0$ of finite index $N$.
Point stabilizers are finite, since they are discrete subgroups
of compact isotropy groups. Each embeds in $\Lambda/\Lambda_0$,
so its order is at most $N$. Proposition~\ref{prop:smooth}
gives a uniform gap for $\Lambda$ and all its subgroups.
This applies to lattices in the isometry group of a symmetric
space of noncompact type: finite generation and Selberg's lemma
give such a $\Lambda_0$, see Witte Morris
\cite[Theorems 4.7.10 and 4.8.2]{Mor15}, after passing to the
identity component and taking a normal core.
Every locally elliptic subgroup is consequently finite by
Corollary~\ref{cor:smooth}, without a finite-generation assumption.
Discreteness is essential here: a cyclic irrational rotation
has no nontrivial finite subgroup but has no uniform elliptic gap.
\end{remark}

\begin{corollary}[Finite image without finite generation]\label{cor:smooth}
Let $C>0$, let $n\geq1$ be an integer, and let a group $G$ act
isometrically and locally elliptically on an $n$-dimensional
Hadamard manifold $M$. Put $\Gamma=\operatorname{im}(G\to\Isom(M))$.
If the action satisfies $\UEG(C)$, then
$|\Gamma|\leq\mathcal B_n(C)$ and $G$ fixes a point of $M$.
The same finiteness and fixed-point conclusions hold if the orders
of the isometries in $\Gamma$ are finite and uniformly bounded.
\end{corollary}
\begin{proof}
Suppose first that the action satisfies $\UEG(C)$, and put
$B=\mathcal B_n(C)$. Every finitely generated subgroup of $\Gamma$
fixes a point by Theorem~\ref{thm:main} and has at most $B$ elements
by Proposition~\ref{prop:smooth}. If $\Gamma$ had a finite subset
of cardinality greater than $B$, the subgroup generated by that
subset would contradict this bound. Thus $\Gamma$ is finite and
fixes the circumcenter of any orbit. The bounded-order conclusion
follows from Proposition~\ref{prop:smooth}.
\end{proof}

\begin{proposition}[Polynomial stabilizer bound]\label{prop:jordan}
Let $n\geq1$, let $M$ be an $n$-dimensional Hadamard manifold, and
let $\Gamma\leq\Isom(M)$ satisfy $\UEG(C)$ with $C\geq1/2$.
With $J(n)$ and $Q(C)$ as in Lemma~\ref{lem:orthogonal},
\[
 |\Gamma_p|\leq J(n)
 \max_{0\leq r\leq\lfloor n/2\rfloor}2^{n-2r}Q(C)^r
 \qquad(p\in M).
\]
If the action is locally elliptic, the same bound holds for
$|\Gamma|$.
\end{proposition}
\begin{proof}
The first assertion is the polynomial term of
\eqref{eq:smooth-bound}, obtained from Lemma~\ref{lem:orthogonal}.
For a locally elliptic action, Corollary~\ref{cor:smooth}
gives a common fixed point $p$, so $\Gamma=\Gamma_p$.
\end{proof}

\begin{remark}[Size of the smooth estimates]\label{rem:smooth-sharpness}
The bounded-order constant in Proposition~\ref{prop:smooth}
is sharp: the cyclic rotation group of
order $m\geq 2$ on $\R^2$ has optimal gap constant
$(2\sin(\pi/m))^{-1}$. For fixed $n$,
Proposition~\ref{prop:jordan} gives
$|\Gamma_p|=O_n(C^{\lfloor n/2\rfloor})$.
This exponent is optimal. With
$r=\lfloor n/2\rfloor$, the group
$(\mathbb{Z}/m\mathbb{Z})^r$ acting by independent rotations
on $r$ orthogonal planes has order $m^r$
and, for $r\geq1$, optimal gap constant
$(2\sin(\pi/m))^{-1}\sim m/(2\pi)$.
For $n=1$, a reflection group realizes
the constant order of growth. No optimal
multiplicative constant is asserted.
\end{remark}

{\bf Acknowledgements.}
The author is funded by the New Zealand Marsden Fund through grant
UOA-2122. OpenAI's Astra 6.0 was used to assist in reviewing and revising mathematical arguments and their presentation, and in preparing the introduction and bibliography. All AI-assisted content was verified and revised by the author.

\end{document}